\documentclass[a4paper]{article}
\usepackage[latin1]{inputenc}
\usepackage[T1]{fontenc}
\usepackage[francais,english]{babel}
\usepackage{amsmath}
\usepackage{amssymb}
\usepackage{theorem}
\newtheorem{thm}{Theorem}
\newtheorem{Prop}{Proposition}
\newtheorem{Def}{Definition}
\newtheorem{Lem}{Lemma}
\newtheorem{Cor}{Corollary}
\newtheorem{rk}{Remark}
\newcommand{\R}{\mathbb{R}}
\newcommand{\C}{\mathbb{C}}

\author{
Karine \textsc{Beauchard}\footnote{CMLA, ENS Cachan, CNRS, Universud,
61 avenue du Pr\'{e}sident Wilson, F-94230 Cachan, France,
email: Karine.Beauchard@cmla.ens-cachan.fr}
Jean-Michel \textsc{Coron}\footnote{Institut universitaire de France and Universit\'{e}
Pierre et Marie Curie-Paris 6, UMR 7598 Laboratoire Jacques-Louis
Lions, Paris, F-75005 France,
email: coron@ann.jussieu.fr}
Pierre \textsc{Rouchon}\footnote{Mines ParisTech,
Centre Automatique et Syst\`{e}mes,
Math\'{e}matiques et Syst\`{e}mes,
60, boulevard Saint-Michel,
75272 Paris CEDEX, France,
email: Pierre.Rouchon@mines-paristech.fr}
\thanks{KB, JMC and PR were partially supported by the ``Agence Nationale de la Recherche'' (ANR),
Projet Blanc C-QUID number BLAN-3-139579}
}

\title{Controllability issues  for continuous-spectrum systems and
ensemble controllability of Bloch Equations}
\date{}

\begin{document}
\maketitle

\begin{abstract}

We study the controllability of the Bloch equation, for an ensemble of
non interacting half-spins, in a static magnetic field, with dispersion in the
Larmor frequency. This system may be seen as a prototype for infinite dimensional
bilinear systems with continuous spectrum, whose controllability is not well
understood. We provide several mathematical answers, with discrimination between
approximate and exact controllability, and between finite time or
infinite time controllability: this system is not exactly controllable
in finite time $T$ with bounded controls in $L^2(0,T)$,
but it is approximately controllable in $L^\infty$ in finite time
with unbounded controls in $L^{\infty}_{loc}([0,+\infty))$.
Moreover, we propose explicit controls realizing the asymptotic
exact controllability to a uniform state of spin $+1/2$ or $-1/2$.

\end{abstract}



\section{Introduction}
\label{section:Intro}

\subsection{Studied system, bibliography}

Most controllability results available for infinite dimensional systems
are related to systems with discrete spectra.
As far as we know, very few  controllability studies consider
systems admitting a continuous part in their spectra.
In~\cite{mirrahimi-IHP09} an approximate controllability result is given
for a system with mixed  discrete/continuous spectrum:
the Schr\"{o}dinger  partial differential equation of a quantum particle
in an N-dimensional decaying potential  is shown to be approximately controllable
(in infinite time) to the ground bounded state when the  initial state is a
linear superposition of bounded states.

In~\cite{Li-Khaneja:cdc06,Li-Khaneja:PRA06,Li-Khaneja:ieee09} a controllability
notion, called ensemble controllability, is  introduced and  discussed for quantum
systems  described by a family of ordinary differential equations (Bloch equations)
depending continuously on a finite number of scalar parameters and with a finite
number of control inputs. Ensemble controllability means that it is possible
to find open-loop controls that compensate for the dispersion in these scalar
parameters: the goal  is to simultaneously steer a continuum of systems between
states of interest with the same control input.
 The articles \cite{Li-Khaneja:cdc06,Li-Khaneja:PRA06,Li-Khaneja:ieee09}  highlight, for three common dispersions in NMR
spectroscopy, the role of Lie algebras and non-commutativity in the design
of a compensating control sequence and  consequently  in the  characterization of
ensemble controllability.

Such continuous family of  ordinary differential systems sharing the same
control inputs  can be seen as the prototype of infinite dimensional
systems with purely continuous spectra. The goal of this paper is to
show that  the very interesting controllability analysis of~\cite{Li-Khaneja:cdc06,Li-Khaneja:PRA06,Li-Khaneja:ieee09}
can be completed  by  functional analysis methods developed for infinite
dimensional systems governed by partial differential equations
(see, e.g., \cite{coron:book} for samples of these methods).

We focus here  on  one of the three dispersions cases treated in~\cite{Li-Khaneja:cdc06,Li-Khaneja:PRA06,Li-Khaneja:ieee09}.
We consider an ensemble of non interacting half-spins
in a static field {\scriptsize $\left(
                     \begin{array}{c}
                       0 \\
                       0 \\
                       B_0 \\
                     \end{array}
                   \right)
$} in $\mathbb{R}^3$,  subject to a transverse radio frequency field {\scriptsize $\left(
                                                                    \begin{array}{c}
                                                                      v(t) \\
                                                                      -u(t) \\
                                                                      0 \\
                                                                    \end{array}
                                                                  \right)
$} in $\mathbb{R}^3$ (the control input). The ensemble of half-spins is
described  by the magnetization vector $M\in\mathbb{R}^3$ depending on time $t$ but also on the Larmor frequency $\omega=-\gamma B_0$ ($\gamma$ is the gyromagnetic
ratio). It obeys to  the Bloch equation:
\begin{equation} \label{eq:syst-Bloch}
\begin{array}{l}
\displaystyle  \frac{\partial M}{\partial t}(t,\omega) =
\left( \begin{array}{ccc}
0 & - \omega & v(t) \\
\omega & 0   & -u(t) \\
-v(t)   & u(t)   &  0
\end{array} \right)
M(t,\omega),\quad
(t,\omega) \in [0,+\infty) \times (\omega_{*},\omega^{*}),
\end{array}
\end{equation}
where $- \infty \leqslant \omega_{*} < \omega^{*} \leqslant + \infty$ are given .
With the notations
\begin{equation} \label{def-Omega-xyz}
\begin{array}{ccc}
\Omega_x := \left(
\begin{array}{ccc}
0 & 0 & 0 \\ 0 & 0 & -1 \\ 0 & 1 & 0
\end{array}
\right),
&
\Omega_y := \left(
\begin{array}{ccc}
0 & 0 & 1 \\ 0 & 0 & 0 \\ -1 & 0 & 0
\end{array}
\right),
&
\Omega_z := \left(
\begin{array}{ccc}
0 & -1 & 0 \\ 1 & 0 & 0 \\ 0 & 0 & 0
\end{array}
\right),
\end{array}
\end{equation}
the system (\ref{eq:syst-Bloch}) can be written
\begin{equation} \label{Bloch-ABC}
\displaystyle \frac{\partial M}{\partial t}(t,\omega)=
\left(\omega \Omega_z + u(t) \Omega_x + v(t) \Omega_y \right) M(t,\omega),
\quad
(t,\omega) \in [0,+\infty) \times (\omega_{*},\omega^{*}).
\end{equation}
It is a bilinear control system in which, at time $t$,
\begin{itemize}
\item the state is $\left(M(t,\omega)\right)_{\omega\in(\omega_{*},\omega^{*})}$; for each $\omega$, $M(t,\omega)\in  \mathbb{S}^{2}$, the unit sphere of $\R^3$,
\item the two control inputs  $u(t)$ and $v(t)$ are real.
\end{itemize}
In the sequel, we denote by $e_k$, the $\R^3$-vector of coordinates $(\delta_{ki})_{i\in\{1,2,3\}}$.
Thus, we study the simultaneous controllability of a continuum
of ordinary differential equations, with respect to a parameter
$\omega$ that belongs to an interval $(\omega_*,\omega^*)$.
Notice that, when $v=u=0$, the spectrum of this system is made
by the union of the two  segments, $i (\omega_*,\omega^*)$
and $-i (\omega_*,\omega^*)$, belonging to the imaginary axis.\par

The pioneer articles \cite{Li-Khaneja:cdc06,Li-Khaneja:PRA06,Li-Khaneja:ieee09} provide convincing arguments indicating why
the system (\ref{Bloch-ABC}) is ensemble  controllable
(i.e. approximately controllable in $L^2((\omega_*,\omega^*),\mathbb{S}^2)$)
with unbounded and also bounded controls, when $\omega_*$ and $\omega^*$
are finite. Here, we provide several mathematical  results that complete these ensemble
controllability results with discriminations between  approximate or exact
controllability  and between finite or infinite time (asymptotically) controllability.

\subsection{Controllability issues}

Let us recall a famous non controllability result for infinite dimensional
bilinear systems due to Ball, Marsden and Slemrod \cite{BMS}. This result
concerns general systems of the form
\begin{equation} \label{blsyst}
\frac{dw}{dt} = \mathcal{A} w + p(t) \mathcal{B}w
\end{equation}
where the state is $w$ and the control is $p:[0,T] \rightarrow \mathbb{R}$.

\begin{thm} \label{thm-BMS}
Let $X$ be a Banach space with $\text{dim}(X)=+\infty$,
$\mathcal{A}$ generate a $C^0$-semigroup of bounded operators on $X$ and
$\mathcal{B}:X \rightarrow X$ be a bounded operator. For $w_0 \in X$,
$w(t;p,w_0)$ denotes the unique solution of (\ref{blsyst})
with $p \in L^1_{loc}([0,+\infty))$ and $w(0)=w_0$. The reachable set from $w_0$
$$\mathcal{R}(w_0) := \{ w(t;p,w_0); t \geqslant 0, p \in L^r_{loc}([0,+\infty)),
r>1\}$$
is contained in a countable union of compact subsets of $X$ and, in particular,
it has an empty interior in $X$. Thus (\ref{blsyst}) is not controllable in $X$
with controls in $\cup_{r>1} L^r_{loc}([0,+\infty))$.
\end{thm}
We cannot apply directly here this result since the spaces $X=L^2((\omega_*,\omega^*),\mathbb{S}^2)$ or $C^0([\omega_*,\omega^*],\mathbb{S}^2)$ where the Cauchy problem is well-defined  are not vector spaces.
In order to get an interesting result for the Bloch equation,
one  needs extensions of the above result to Banach manifolds. (This has been done in~\cite{Turinici} when the manifold is the unit sphere of a Hilbert space.)   For~\eqref{def-Omega-xyz}, the situation is similar to the one described in Theorem \ref{thm-BMS}. In Theorem~\ref{analytic}, we show that for any analytic initial condition $M_0(\omega)$, the reachable set in finite time $T>0$
from $M_0$ with controls in $L^2(0,T)$ only contains analytic functions
of $\omega$. Thus, the reachable set (from an analytic initial dara) has an empty interior in
$L^2((\omega_*,\omega^*),\mathbb{S}^2)$, which is a natural space for the Cauchy problem.

However, for~\eqref{def-Omega-xyz}, the obstruction to exact controllability given by Theorem~\ref{analytic} has much stronger consequence than the obstruction described by Theorem \ref{thm-BMS} which is, in fact,  a rather weak non controllability
result. Indeed, it does not prevent the reachable set from being dense in $X$
(approximate controllability in $X$). For example, this is the case for the 1D beam equation
$$\left\lbrace \begin{array}{l}
u_{tt}+u_{xxxx}+p(t) u_{xx}=0, \quad x \in (0,1),~ t \in (0,+\infty), \\
u=u_x=0 \text{ at } x=0,1,
\end{array}\right.$$
in which the state is $(u,u_t)$ and the control is $p$.
Theorem \ref{thm-BMS} ensures that this system is not exactly controllable in
$H^2_0 \times L^2 (0,1)$ with controls in $L^r_{loc}([0,+\infty))$, $r>1$.
However, it is proved in \cite{KB-beam} that this system is exactly
controllable in $H^{5+} \times H^{3+} (0,1)$ with controls in $H^1_0(0,T)$,
at least locally around a stationnary trajectory.
Similarly, Turinici's generalization~\cite{Turinici} of Theorem \ref{thm-BMS} applies to
1D Schr\"{o}dinger equations of the form
$$\left\lbrace\begin{array}{l}
i \frac{\partial \psi}{\partial t}= - \frac{\partial^2 \psi}{\partial x^2}
- u(t) \mu(x) \psi,\quad  x \in (0,1),~ t \in (0,+\infty), \\
\psi(t,0)=\psi(t,1)=0
\end{array}\right.$$
where the state is $\psi$, the control is $u$ and $\mu \in C^\infty([0,1])$.
It proves that this system is not exactly controllable in $H^2((0,1),\mathbb{C})$
with controls in $L^2_{loc}([0,+\infty))$. However it is proved in \cite{05Beauchard,2006-jfa} that this system, with $\mu(x)=(x-1/2)$ is exactly controllable
in $H^7((0,1),\mathbb{C})$ with controls in $H^1_0(0,T)$,
locally around the eigenstates, for $T$ large enough.

The conclusion of~\cite{05Beauchard,KB-beam,2006-jfa} is that, sometimes, the negative result of Theorem~\ref{thm-BMS} is only due to a bad choice of functional spaces that do not
allow the controllability; but positive controllability results may
be expected in different functional spaces.
Therefore, one may still hope to prove the exact controllability of the Bloch equation
in some well chosen functional spaces. We will see in this article that
it is not the case: the Bloch equation is not exactly controllable
in a much stronger sense than the one of Theorem \ref{thm-BMS}.

Indeed, we will prove that, when $(\omega_*,\omega^*)=(-\infty,+\infty)$,
the reachable set  (in finite time and with small controls) from $M_0 \equiv e_3$ is a \textbf{submanifold} of some functional space,
that does not coincide with one of its tangent spaces.
When the domain $(\omega_*,\omega^*)$ is a bounded interval of $\mathbb{R}$,
we will see that there exist \textbf{analytic} targets, arbitrarily close to $e_3$
that cannot be reached exactly from $e_3$ with bounded controls in $L^2(0,T)$.
Thus, the non controllability of (\ref{Bloch-ABC}) is not related to
a regularity problem and this equation corresponds to a very different situation
from~\cite{05Beauchard,KB-beam,2006-jfa}.

\subsection{Outline and open problems}

In section \ref{Section:Linearise}, we study the linearized system of (\ref{Bloch-ABC})
around the steady-state   $(M \equiv e_3, (u,v) \equiv 0)$
with $-\infty <\omega_*<\omega^*<+\infty$.
This system is  shown to be approximately controllable in
$C^{0}([\omega_*,\omega^*],\mathbb{R}^3)$,
in any finite time $T$, with unbounded controls $(u,v) \in C^{\infty}_{c}((0,T), \R^2)$.
But it is not exactly controllable
neither in finite time nor in infinite time.
Moreover, for any reachable target, there exists only one control which steers
the control system to the target.

In section \ref{section:R-T}, we study the exact controllability
of the nonlinear system (\ref{Bloch-ABC}),
locally around $M \equiv e_3$, in finite time.
 First, we prove that the simultaneous exact controllability with respect
to $\omega$ in the whole space $\mathbb{R}$ (i.e. $\omega_*=-\infty$,
$\omega^*=+\infty$) does not hold with bounded controls.
Indeed, for every time $T>0$, the reachable set from $M_0 \equiv e_3$
with bounded controls in
$L^2(0,T)$ is a strict submanifold (of some functional space) that does not coincide
with one of its tangent space.
Then, with an analyticity argument, we deduce that the simultaneous
exact controllability with respect to $\omega$ in a bounded interval
$(\omega_*,\omega^*)$, $-\infty<\omega_*<\omega^*<+\infty$,
does not hold neither.

The exact controllability of (\ref{Bloch-ABC}) being impossible
with bounded controls, in sections
\ref{Li-Khaneja} and \ref{Section:asymptotic}, we investigate the exact
controllability of (\ref{Bloch-ABC}) with unbounded controls.

In section \ref{Li-Khaneja}, completing the arguments of
~\cite{Li-Khaneja:cdc06,Li-Khaneja:PRA06,Li-Khaneja:ieee09}, we prove the ensemble controllability
of (\ref{Bloch-ABC}): any measurable initial condition
$M_0:(\omega_*,\omega^*) \rightarrow \mathbb{S}^2$
can be steered approximately in $L^2(\omega_*,\omega^*)$ to $e_3$.
This approximate controllability indeed holds for stronger norms,
for instance $\|.\|_{L^{\infty}}$ and $\| . \|_{H^s}$, $\forall s \in (0,1)$.
The controls used to realize this
motion are sequences of pulses presented in ~\cite{Li-Khaneja:cdc06}
(but one may also use controls in $L^{\infty}_{loc}([0,+\infty))$)
and the proof relies on non-commutativity and functional analysis.

In section \ref{Section:asymptotic}, we propose other explicit
unbounded controls realizing the asymptotic local (exact)
controllability to $e_3$, simultaneously with respect to $\omega$
in a bounded interval. Here, the proof relies on Fourier analysis.

 Finally, in section \ref{sec:Comparaison}  ,
we compare the feasibility, the time and the cost of the two controllability processes
presented in sections \ref{Li-Khaneja} and \ref{Section:asymptotic},
on particular motions.

Let us emphasize that the behavior of the nonlinear system around $e_3$ is
very different from the one of the linearized system around $e_3$. Indeed,
\begin{itemize}
\item first, the linearized system is not asymptotically zero controllable
whereas the nonlinear system is asymptotically locally controllable to $e_3$,
\item then, as seen in section \ref{Section:Linearise},
for the linearized system and for any reachable target,
only a single  control works, whereas for the nonlinear system and
for any initial condition, many controls
allow to reach exactly $e_3$ (in infinite time).
\end{itemize}
Thus, the nonlinearity allows to recover controllability.
 Finally, let us mention some open problems.

In section \ref{section:R-T}, we prove the non exact controllability to $e_3$ with
bounded controls, in finite time, because the reachable set is a submanifold.
The equation of this submanifold and the validity of the same negative
result in infinite time (i.e. the non asymptotic exact controllability to $e_3$
with bounded controls) are open problems.

In section \ref{Section:asymptotic}, we prove the exact controllability to $e_3$
with unbounded controls, in infinite time. The validity of the same result
in finite time is also open.

Before starting the mathematical study
let us introduce some notations that will be used in all the paper.
We write
\begin{equation} \label{Def:x,y,z}
M(t,\omega):=\left( \begin{array}{c}
x(t,\omega) \\ y(t,\omega)\\ z(t,\omega)
\end{array}\right),
\end{equation}
\begin{equation} \label{Def:Z,w}
\begin{array}{ll}
Z(t,\omega):=(x+iy)(t,\omega),
&
w(t):=(-v+iu)(t).
\end{array}
\end{equation}
Thus, when, for some time $T>0$, $z(t,\omega)>0$ on $(0,T) \times (\omega_*,\omega^*)$,
then the system (\ref{Bloch-ABC}) implies
\begin{equation} \label{eq:Z}
\frac{\partial Z}{\partial t}(t,\omega)=
i \omega Z(t,\omega) - w(t) \sqrt{1-|Z(t,\omega)|^2},\quad
(t,\omega) \in (0,T) \times (\omega_*,\omega^*),
\end{equation}
so
$$Z(t,\omega)=\left(
Z_0(\omega) - \int_{0}^{t} w(\tau) \sqrt{1-|Z(\tau,\omega)|^2} e^{-i\omega \tau}
d\tau \right) e^{i \omega t},\quad
(t,\omega) \in (0,T) \times (\omega_*,\omega^*).$$
Unless otherwise specified, the functions considered are complex valued and, for example, we write $L^2(\R)$ for $L^2(\R,\C)$. When the functions considered are real valued we specify it and, for example, we write $L^2(\R,\R)$.

\section{Linearized system around $(M \equiv e_3, (u,v) \equiv 0)$}
\label{Section:Linearise}

In this section $-\infty < \omega_* < \omega^* < + \infty$.
We are interested in the linearized system of (\ref{Bloch-ABC})
around $(M \equiv e_3, (u,v) \equiv 0)$, or, equivalently
in the linearized system of (\ref{eq:Z}) around $(Z \equiv 0, w \equiv 0)$,
\begin{equation} \label{linearise}
\dot{Z}(t,\omega)=i \omega Z(t,\omega) - w(t), \quad Z(0,\omega)=Z_0(\omega),
\end{equation}
whose solution is
\begin{equation} \label{lin-sol}
Z(t,\omega)= \Big( Z_0(\omega) -
\int_{0}^{t} w(\tau) e^{-i\omega \tau} d\tau \Big) e^{i \omega t}.
\end{equation}
We prove its  non exact controllability and its
approximate controllability with unbounded controls.

\subsection{Non exact controllability}

We denote by $\mathcal{F}$ the 1-D Fourier transform:
$$
\mathcal{F}(w)(\omega)=\int_\R w(t)e^{-i\omega t}dt.
$$
When a function is defined on $I\subset \R$, we extend it by $0$ on $\R\setminus I$.
One has the following proposition.
\begin{Prop} \label{Linearise-neg}
Let $T \in (0,+\infty)$. The reachable set from $Z_0=0$
for (\ref{linearise}) with controls $w \in L^1(0,T)$ is
$$\{ Z(T) ;\, w \in L^1(0,T) \} =
\mathcal{F}[ L^1(-T,0)].$$

The set of initial conditions $Z_0$ that are asymptotically zero
controllable with controls $w \in L^1(0,+\infty)$ for
(\ref{linearise}) is $\mathcal{F}[L^1(0,+\infty)]$.

 For every $Z_0 \in \mathcal{F}[L^1(0,+\infty)]$,
the function $w:= \mathcal{F}^{-1}[Z_0]$ is the unique control
in $L^1(0,+\infty)$ that steers the control system (\ref{linearise}) from $Z_0$ to $0$.
\end{Prop}

\textbf{Proof of Proposition \ref{Linearise-neg}: }
The two first statements are direct consequences of the explicit expression
(\ref{lin-sol}). Concerning the third statement, it is sufficient to prove
that if $w \in L^1(0,+\infty)$ and if $\mathcal{F}[w] \equiv 0$
on $(\omega_*,\omega^*)$, then $w=0$. Let $w$ be such a function and consider
$\varphi: \mathbb{C}_{+} \cup \mathbb{C}_{-} \cup (\omega_*,\omega^*)
\rightarrow \mathbb{C}$, defined by
$$\varphi(z):=
\left\lbrace  \begin{array}{l}
\mathcal{F}[w](z) , \text{ if } z \in \mathbb{C}_{-}, \\
\overline{\mathcal{F}[w](\overline{z})}, \text{ if } z \in \mathbb{C}_{+},\\
0 , \text{ if } z \in (\omega_*,\omega^*),
\end{array}\right.$$
where $\mathbb{C}_{+} := \{ z \in \mathbb{C} ;\, \Im(z)>0 \}$
and $\mathbb{C}_{-} := \{ z \in \mathbb{C} ;\, \Im(z)<0 \}$.
Then $\varphi$ is holomorphic on $\mathbb{C}_{+}$ and on $\mathbb{C}_{-}$ and
continuous on $\mathbb{C}_{+} \cup \mathbb{C}_{-} \cup (\omega_*,\omega^*)$,
so it is holomorphic on $\mathbb{C}_{+} \cup \mathbb{C}_{-} \cup (\omega_*,\omega^*)$.
Since $\varphi$ vanishes on $(\omega_*,\omega^*)$, then $\varphi \equiv 0$.
Thus $w=0$. $\Box$

\subsection{Approximate controllability with unbounded controls}

\begin{Prop} \label{Weierstrass}
Let $T>0$, $Z_f \in C^{0}([\omega_*,\omega^*])$ and $\eta >0$.
There exists $w \in C^{\infty}_{c}((0,T))$ such that the solution of
(\ref{linearise}) with $Z_0=0$
satisfies $\|Z(T)-Z_f\|_{L^{\infty}(\omega_*,\omega^*)} < \eta$.
\end{Prop}

\textbf{Proof of Proposition \ref{Weierstrass}:}
Let $T>0$, $Z_f \in C^{0}([\omega_*,\omega^*])$ and $\eta >0$.
Thanks to the Weierstrass theorem, there exists a polynomial $P \in \mathbb{C}[X]$ such that
$$\| Z_{f}(\omega) e^{-i\omega \frac{T}{2}} -
P(i \omega) \|_{L^{\infty}(\omega_*,\omega^*)} < \frac{\eta}{2}.$$
Applying the control $w(t):= - P(\partial_{t}) \delta_{T/2}(t)$ in
(\ref{lin-sol}) with $Z_0=0$, we get
$$Z(T,\omega) = - \mathcal{F}[w](\omega) e^{i \omega T}
= P(i \omega) e^{-i\omega \frac{T}{2}} e^{i \omega T}
= P(i \omega) e^{i\omega \frac{T}{2}}$$
thus,
$$\| Z(T) - Z_{f} \|_{L^{\infty}(\omega_*,\omega^*)} < \frac{\eta}{2}.$$
Now, let us smooth this control candidate in order to provide a smooth control.
Let $g \in C^{\infty}_{c}((-1,1),\mathbb{R}_+)$ such that $\int_{\mathbb{R}} g =1$.
 For $\epsilon \in (0,T/2)$, the function
$$g_{\epsilon}(t) := \frac{1}{\epsilon} g \left( \frac{t-T/2}{\epsilon} \right)$$
is supported in $(0,T)$. Applying the control
$w_{\epsilon}(t):=  - P(\partial_{t}) g_{\epsilon}(t)$ in (\ref{linearise}), we get
$$Z(T,\omega)  = P(i \omega) \hat{g_{\epsilon}} (\omega) e^{i \omega T}.$$
Noticing that
$$\hat{g_{\epsilon}}(\omega) - e^{-i \omega \frac{T}{2}}
= \int_{-1}^{1} g(y) [ e^{-i \omega \epsilon y} -1 ] dy~ e^{-i\omega \frac{T}{2}},$$
we get
$$\| P(i \omega) [ \hat{g_{\epsilon}} (\omega) - e^{-i \omega \frac{T}{2}} ]
\|_{L^{\infty}(\omega_*,\omega^*)}
\leqslant
\| P(i\omega) \|_{L^{\infty}(\omega_*,\omega^*)}
\| \hat{g_{\epsilon}} (\omega) - e^{-i \omega \frac{T}{2}} \|_{L^{\infty}(\omega_*,\omega^*)}
\rightarrow 0$$
when $\epsilon \rightarrow 0$. Thus, for $\epsilon$ small enough,
$$\| Z(T) - Z_{f} \|_{L^{\infty}(\omega_*,\omega^*)} < \eta. \Box$$

\section{Non exact controllability with bounded controls}
\label{section:R-T}

In this section, we study the reachable set from $M(0) \equiv e_{3}$
for (\ref{eq:syst-Bloch})
with bounded controls $(u,v) \in L^{2}((0,T),\mathbb{R}^{2})$.
Notice that, when $M(0) \equiv e_{3}$ and
$w$ is small enough in $L^{1}(0,T)$,
then $z(t,\omega)>0$ for every $(t,\omega) \in (0,T) \times \mathbb{R}$ and
\begin{equation} \label{eq:Z-sol-faible}
Z(t,\omega)=-\int_{0}^{t} w(\tau) \sqrt{1-|Z(\tau,\omega)|^{2}}
e^{-i\omega \tau} d\tau e^{i \omega t},
\forall (t,\omega) \in [0,T] \times \mathbb{R}.
\end{equation}

\subsection{Case $\omega_*=-\infty$, $\omega^*=+\infty$}

In this section, we take $\omega_*=-\infty$, $\omega^*=+\infty$.
In a first subsection we precise the functional framework in which
(\ref{eq:Z-sol-faible}) is well posed. In a second subsection,
we prove that the reachable set from zero, with bounded controls $(u,v)$
in $L^2((0,T),\mathbb{R}^2)$ is a strict submanifold
of some functional space that does not coincide with its tangent
space at zero. In particular,  (\ref{Bloch-ABC}) is not locally controllable
with bounded controls $(u,v)$ in $L^2((0,T),\mathbb{R}^2)$.

\subsubsection{Solutions on $[0,T]$}

\begin{Prop} \label{Prop:existence-sol-faible}
Let $T>0$ and $R := 1/(2\sqrt{T})$.
 For every $w \in L^{2}(0,T)$ with $\|w\|_{L^{2}(0,T)} < R$,
there exists a unique $Z \in C^{0}([0,T],L^{2}(\mathbb{R}))
\cap C^{0}_{b}([0,T]\times \mathbb{R})$ solution of (\ref{eq:Z-sol-faible})
and it satisfies
\begin{equation} \label{Z-borne-L-infini}
\| Z \|_{L^{\infty}((0,T) \times \mathbb{R})} \leqslant \sqrt{T} \|w\|_{L^{2}(0,T)},
\end{equation}
\begin{equation} \label{Z-borne-L-2}
\| Z \|_{C^{0}([0,T],L^{2}(\mathbb{R}))} \leqslant 2 \sqrt{2\pi}\|w\|_{L^{2}(0,T)}.
\end{equation}
Moreover, for every $w_1, w_2 \in L^{2}(0,T)$ with
$\|w_1\|_{L^{2}(0,T)}< R$ and $\|w_2\|_{L^{2}(0,T)} < R$, we have
\begin{equation} \label{continuity-L-infini}
\| Z_1 - Z_2 \|_{L^{\infty}((0,T) \times \mathbb{R})}
\leqslant
2\sqrt{T} \| w_1 - w_2 \|_{L^{2}(0,T)},
\end{equation}
\begin{equation} \label{continuity-L-2}
\| Z_1 - Z_2 \|_{C^{0}([0,T],L^{2}(\mathbb{R}))}
\leqslant
4 \sqrt{2\pi}\| w_1 - w_2 \|_{L^{2}(0,T)},
\end{equation}
where, for $j\in \{1,2\}$, $Z_j$ denotes the unique solution of (\ref{eq:Z-sol-faible}) for $w:=w_j$.
\end{Prop}

\noindent \textbf{Proof of Proposition \ref{Prop:existence-sol-faible}:}
Let $T>0$ and $c:=1/\sqrt{3}$, which is chosen so that
\begin{equation} \label{def:c}
|f'(x)|  \leqslant c , \forall x \in [0,1/2],
\text{  where  } f(x):=\sqrt{1-x^{2}}.
\end{equation}
Let  $w \in L^{2}(0,T)$ be such that $\|w\|_{L^{2}(0,T)} < R$.
We apply the Banach fixed point theorem to the map $\Theta$ defined on
$$B :=
C^{0}([0,T],L^{2}(\mathbb{R})) \cap
C^{0}([0,T] \times \mathbb{R},B_{\mathbb{C}}(0,1/2))$$
by $\Theta(\xi)=Z$ where
\begin{gather*}
B_{\mathbb{C}}(0,1/2)):=\{\xi \in \C; \, |z|\leqslant 1/2\},
\\
Z(t,\omega)=
-\int_{0}^{t} w(\tau) \sqrt{1-|\xi(\tau,\omega)|^{2}}
e^{-i\omega \tau} d\tau e^{i \omega t},
\forall (t,\omega) \in [0,T] \times \mathbb{R}.
\end{gather*}
Note that $B$ is a nonempty closed subset of the Banach space $C^0([0,T], L^2(\R))$.

\emph{First step: $\Theta$ takes its values in $B$ because $R \sqrt{T} \leqslant 1/2$.}
\\
Let $\xi \in B$ and $Z:=\Theta(\xi)$. The Cauchy-Schwarz inequality leads to
\begin{equation} \label{Z-xi-borne-L-infini}
|Z(t,\omega)| \leqslant \|w\|_{L^{1}(0,T)}
\leqslant \sqrt{T} \|w\|_{L^{2}(0,T)}
\leqslant \sqrt{T} R
\leqslant 1/2.
\end{equation}
Thus $Z \in C^{0}([0,T]\times \mathbb{R},B_{\mathbb{C}}(0,1/2))$.
Thanks to the decomposition
\begin{equation} \label{Z-decomposition}
Z(t,\omega)= - \mathcal{F}[ \tau_{-t} w_{|[0,t]} ](\omega)
-\int_{0}^{t} w(\tau) \Big( \sqrt{1-|\xi(\tau,\omega)|^{2}} - 1 \Big)
e^{-i\omega \tau} d\tau e^{i \omega t},
\end{equation}
where $\tau_{a}(\varphi)(s):=\varphi (s-a)$,
thanks to the Plancherel theorem, (\ref{def:c}) and the Cauchy-Schwarz inequality, we get
\begin{equation} \label{calcul-type-L2}
\begin{array}{ll}
\| Z(t) \|_{L^{2}(\mathbb{R})}
& \leqslant \sqrt{2\pi} \| w \|_{L^{2}(0,T)} +
\left( \int_{\mathbb{R}}
\Big| \int_{0}^{t} w(\tau) \Big( \sqrt{1-|\xi(\tau,\omega)|^{2}} - 1 \Big)
e^{-i\omega \tau} d\tau \Big|^{2} d\omega \right)^{1/2}
\\
& \leqslant  \sqrt{2\pi}\| w \|_{L^{2}(0,T)} +
\left( \int_{\mathbb{R}} \|w\|_{L^{2}(0,T)}^{2} \int_{0}^{t}
c^{2} |\xi(\tau,\omega)|^{2} d\tau d\omega
\right)^{1/2}
\\
& \leqslant \| w \|_{L^{2}(0,T)} \Big(
 \sqrt{2\pi} + c \sqrt{T} \| \xi \|_{C^{0}([0,T],L^{2}(\mathbb{R}))} \Big),
\end{array}
\end{equation}
so $Z(t) \in L^{2}(\mathbb{R})$ for every $t \in [0,T]$.
In the right-hand side of (\ref{Z-decomposition}), the first term belongs to
$C^{0}([0,T],L^{2}(\mathbb{R}))$ and the second term also, as one can prove by
applying the dominated convergence theorem with the following domination,
that holds for every $(t,\omega) \in [0,T] \times \mathbb{R}$,
$$\Big| \int_{0}^{t} w(\tau) \Big( \sqrt{1-|\xi(\tau,\omega)|^{2}} - 1 \Big)
e^{-i\omega \tau} d\tau \Big|
\leqslant
c \|w\|_{L^{2}(0,T)} \left( \int_{0}^{T} |\xi(\tau,\omega)|^{2} \right)^{1/2}.$$

\emph{Second step: $\Theta$ is a contraction on $B$ because $c \sqrt{T} R < 1$.}
\\
Let $\xi_1, \xi_2 \in B$, $Z_1:=\Theta(\xi_1)$ and $Z_2:=\Theta(\xi_2)$.
We have
$$(Z_1-Z_2)(t,\omega)=
- \int_{0}^{t} w(\tau) \Big(
\sqrt{1-|\xi_1(\tau,\omega)|^{2}} - \sqrt{1-|\xi_2(\tau,\omega)|^2} \Big)
e^{-i\omega \tau} d\tau e^{i \omega t}.$$
Using (\ref{def:c}) and  the Cauchy-Schwarz inequality, we get
$$\| Z_1-Z_2 \|_{L^{\infty}((0,T) \times \mathbb{R})}
\leqslant
\| w \|_{L^{1}(0,T)} c \|\xi_1 - \xi_2 \|_{L^{\infty}((0,T)\times \mathbb{R})}
\leqslant
c \sqrt{T} R \|\xi_1 - \xi_2 \|_{L^{\infty}((0,T)\times \mathbb{R})}.$$
Working as in (\ref{calcul-type-L2}), we also get
$$\|(Z_1 - Z_2)(t,.) \|_{L^{2}(\mathbb{R})}
\leqslant
c \sqrt{T} R \| \xi_1 - \xi_2 \|_{C^{0}([0,T],L^{2}(\mathbb{R}))},
\forall t \in [0,T].$$

\emph{Third step: Proof of (\ref{Z-borne-L-infini}) and (\ref{Z-borne-L-2})
thanks to $c \sqrt{T} R \leqslant 1/2$.}
\\
Since $c \sqrt{T} R \leqslant 1/2$,
the inequalities (\ref{Z-borne-L-infini}) and (\ref{Z-borne-L-2})
are consequences of (\ref{Z-xi-borne-L-infini}) and (\ref{calcul-type-L2}) with $\xi=Z$.
\\

\emph{Fourth step: Proof of (\ref{continuity-L-infini}) and (\ref{continuity-L-2})
thanks to $c \sqrt{T} R \leqslant 1/2$.}
\\
Using the decomposition
$$\begin{array}{ll}
(Z_1-Z_2)(t,\omega)=
&
- \int_0^t (w_1 - w_2)(\tau) \sqrt{1 - |Z_1(\tau,\omega)|^{2}} e^{-i\omega \tau} d\tau
e^{i \omega t}
\\ &
- \int_{0}^{t} w_2(\tau) \Big(
\sqrt{1 - |Z_1(\tau,\omega)|^{2}} - \sqrt{1 - |Z_2(\tau,\omega)|^{2}} \Big)
e^{-i\omega \tau} d\tau e^{i \omega t},
\end{array}$$
 the Cauchy-Schwarz inequality and (\ref{def:c}), we get
$$|(Z_1 - Z_2)(t,\omega)| \leqslant
\sqrt{T} \|w_1 - w_2\|_{L^{2}(0,T)} +
c \sqrt{T} R \|Z_1 - Z_2\|_{L^{\infty}((0,T)\times \mathbb{R})},
\forall (t,\omega) \in [0,T] \times \mathbb{R}$$
which, since $c \sqrt{T} R \leqslant 1/2$, leads to
$$\|Z_1 - Z_2 \|_{L^{\infty}((0,T)\times \mathbb{R})}
\leqslant
2 \sqrt{T} \|w_1 - w_2\|_{L^{2}(0,T)}.$$
Using the decomposition
$$\begin{array}{ll}
(Z_1-Z_2)(t,\omega)=
&
- \mathcal{F}[ \tau_{-t}(w_1-w_2)_{|[0,t]} ](\omega)
\\ &
- \int_0^t (w_1 - w_2)(\tau)
\Big( \sqrt{1 - |Z_1(\tau,\omega)|^{2}} - 1 \Big) e^{-i\omega \tau} d\tau e^{i \omega t}
\\ &
- \int_{0}^{t} w_2(\tau) \Big(
\sqrt{1 - |Z_1(\tau,\omega)|^{2}} - \sqrt{1 - |Z_2(\tau,\omega)|^{2}} \Big)
e^{-i\omega \tau} d\tau e^{i \omega t},
\end{array}$$
and working as in (\ref{calcul-type-L2}), we get
$$\begin{array}{ll}
\| (Z_1-Z_2)(t) \|_{L^{2}(\mathbb{R})} \leqslant \sqrt{2\pi}
&
\| w_1 - w_2 \|_{L^{2}(0,T)} +
\|w_1 - w_2 \|_{L^{2}(0,T)} c \sqrt{T}  \|Z_1\|_{C^{0}([0,T],L^{2}(\mathbb{R}))}
\\ & +
\|w_2\|_{L^{2}(0,T)} c \sqrt{T}  \|Z_1-Z_2\|_{C^{0}([0,T],L^{2}(\mathbb{R}))}.
\end{array}$$
Thus, using (\ref{Z-borne-L-2}) and $c \sqrt{T} R \leqslant 1/2$, we get
$$\| Z_1-Z_2 \|_{C^{0}([0,T],L^{2}(\mathbb{R}))}
\leqslant
2 \sqrt{2\pi} \Big( 1 + 2 c \sqrt{T} R  \Big) \| w_1 - w_2 \|_{L^{2}(0,T)}
\leqslant
4\sqrt{2\pi} \| w_1 - w_2 \|_{L^{2}(0,T)}.$$
This shows the existence the existence part of Proposition \ref{Prop:existence-sol-faible} and the uniqueness if one requires that $Z$ takes its values in $B_{\C}(0,1/2)$. The uniqueness without assuming this last assumption can be easily obtained from the previous study by noticing that this study implies that, if two solutions are equal on $[0,\tau]$ with $\tau \in [0,T)$, then there are equal $[0,\tau']$ for $\tau'>\tau$. $\Box$

\subsubsection{Structure of the reachable set from zero in time $T$}

The goal of this section is the proof of the following result, where $$B_R [L^{2}(0,T)]:=\{w\in L^2(0,T);\, \|w\|_{L^2(0,T)}<R\}.$$

\begin{thm} \label{main-thm}
Let $T>0$ and $R :=1/(4\sqrt{3T})$. The image of the end point map
\begin{equation} \label{def:FT}
\begin{array}{cccl}
 F_T : & B_R [L^{2}(0,T)] & \rightarrow & L^{2} \cap C^{0}_{b}(\mathbb{R})\\
      &            w                 & \mapsto     & Z(T,.)
\text{ where } Z \text{ solves } (\ref{eq:Z-sol-faible}),
\end{array}
\end{equation}
is a strict submanifold of $L^{2} \cap C^{0}_{b}(\mathbb{R})$ of infinite codimension
that does not coincide with its tangent space at zero.
\end{thm}

The proof of Theorem \ref{main-thm} relies on the following results
(see \cite[Theorem 73.E and Corollary 73.45, Chapter 73]{Zeidler}).

\begin{thm} \label{Zeidler-thm}
Let $M$ and $N$ be two $C^{k}$-Banach manifolds with chart space over $\mathbb{R}$ and
$k \in \{\infty\}\cup \mathbb{N}\setminus\{0\}$. Let $F:M \rightarrow N$ be a map of class $C^{k}$.
If $F$ is a $C^{k}$ embedding, then $S:=F(M)$ is a submanifold of $N$ and
in particular a $C^{k}$-Banach manifold.
\end{thm}

\begin{thm} Under the same assumptions as in Theorem \ref{Zeidler-thm},
if $F$ is an injective $C^k$ immersion and if $F$ is closed, then $F$ is
a $C^k$ embedding.
\end{thm}

\noindent \textbf{Proof of Theorem \ref{main-thm}:} We take
$M:=B_R [L^{2}(0,T)]$ and $N: =L^{2} \cap C^{0}_{b}(\mathbb{R})$. They
are both $C^{\infty}$-Banach manifolds as open subsets of
Banach spaces. The continuity of $F_T$ is a consequence of
(\ref{continuity-L-infini}) and (\ref{continuity-L-2}).
With similar manipulations as in the proof of (\ref{continuity-L-infini})
and (\ref{continuity-L-2}), one can prove that
$F_T$ is $C^1$ and $dF_T(w).W=\xi(T,.)$ where $\xi$ is defined by
\begin{equation} \label{def:differentielle}
\begin{array}{c}
\begin{array}{ll}
\xi(t,\omega)= &
- \int_{0}^{t} W(\tau) \sqrt{1-|Z(\tau,\omega)|^{2}} e^{-i\omega \tau} d\tau e^{i \omega t}
\\ &
+ \int_{0}^{t} w(\tau)
\displaystyle
\frac{\Re [ \overline{Z(\tau,\omega)} \xi(\tau,\omega) ]}{\sqrt{1-|Z(\tau,\omega)|^{2}}}
e^{-i\omega \tau} d\tau e^{i \omega t},
\forall (t,\omega) \in [0,T] \times \mathbb{R}.
\end{array}
\end{array}
\end{equation}
We use the same notation $c$ as in the previous proof (see (\ref{def:c})).
\\

\emph{First step: $F_T$ is injective on $B_R [L^{2}(0,T)]$
because $6 c \sqrt{T} R < 1$.}
\\
Let $w_1, w_2 \in B_R [L^{2}(0,T)]$ be such that $F_T(w_1)=F_T(w_2)$.
 From
$$\int_{0}^{T} w_1(t) \sqrt{1-|Z_1(t,\omega)|^2} e^{-i\omega t} dt
= \int_{0}^{T} w_2(t) \sqrt{1-|Z_2(t,\omega)|^2} e^{-i\omega t} dt,$$
we deduce
$$\begin{array}{ll}
\mathcal{F}[w_1-w_2](\omega) =
&
\int_{0}^{T} (w_2-w_1)(t) \Big( \sqrt{1-|Z_2(t,\omega)|^2} - 1 \Big) e^{-i\omega t} dt
\\ &
+ \int_{0}^{T} w_1(t) \Big(  \sqrt{1-|Z_2(t,\omega)|^2} -
\sqrt{1-|Z_1(t,\omega)|^2} \Big) e^{-i\omega t} dt.
\end{array}$$
Considering the $L^{2}(\mathbb{R})$-norm of both sides,
using Plancherel equality and working as in (\ref{calcul-type-L2})
we get
$$\begin{array}{ll}
\sqrt{2\pi} \|w_1-w_2\|_{L^{2}(0,T)} \leqslant
&
\|w_2-w_1\|_{L^{2}(0,T)} c \sqrt{T} \|Z_{2}\|_{C^{0}([0,T],L^{2}(\R))}
\\ &
+ \|w_1\|_{L^{2}(0,T)} c \sqrt{T} \| Z_1 - Z_2 \|_{C^{0}([0,T],L^{2}(\R))}.
\end{array}$$
Using (\ref{Z-borne-L-2}) and (\ref{continuity-L-2}), we deduce
$$\sqrt{2\pi}\|w_1-w_2\|_{L^{2}(0,T)}
\leqslant
6  c \sqrt{2\pi}\sqrt{T} R  \|w_1 - w_2\|_{L^{2}(0,T)},$$
which gives the conclusion, because $6 c \sqrt{T} R < 1$.
\\

\emph{Second step: $F_T$ is an immersion because $6 c \sqrt{T} R < 1$.}
\\
Let $w \in B_{R}[L^{2}(0,T)]$ and $W \in L^{2}(0,T)$
be such that $dF_T(w).W=0$. Thanks to (\ref{def:differentielle}), we have
$$\begin{array}{ll}
\mathcal{F}[W](\omega) =
&
- \int_{0}^{T} W(\tau) \Big(
\sqrt{1-|Z(\tau,\omega)|^{2}}-1 \Big) e^{-i\omega \tau} d\tau
\\ &
+ \int_{0}^{T} w(\tau)
\displaystyle \frac{\Re [ \overline{Z(\tau,\omega)} \xi(\tau,\omega) ]}{\sqrt{1-|Z(\tau,\omega)|^{2}}}
e^{-i\omega \tau} d\tau,
\forall \omega \in  \mathbb{R}.
\end{array}$$
Considering the $L^{2}(\mathbb{R})$-norm of both sides and
working as in (\ref{calcul-type-L2}) we get
$$\sqrt{2\pi} \|W\|_{L^{2}(0,T)} \leqslant
\|W\|_{L^{2}(0,T)} c \sqrt{T} \| Z \|_{C^{0}([0,T],L^{2}(\R))} +
\|w\|_{L^{2}(0,T)} c \sqrt{T} \| \xi \|_{C^{0}([0,T],L^{2}(\R))}$$
Admitting the following inequality,
\begin{equation} \label{diff-cont}
\| \xi \|_{C^{0}([0,T],L^{2}(\mathbb{R}))}
\leqslant 4 \sqrt{2\pi}\|W\|_{L^{2}(0,T)},
\end{equation}
and using (\ref{Z-borne-L-2}), we get
$$\sqrt{2\pi} \|W\|_{L^{2}(0,T)} \leqslant 6 c \sqrt{2\pi}\sqrt{T} R \|W\|_{L^{2}(0,T)},$$
which gives the conclusion because $6c \sqrt{T} R < 1$.
Now, let us prove (\ref{diff-cont}). Using the decomposition
$$\begin{array}{ll}
\xi(t,\omega)= &
- \mathcal{F}[ \tau_{-t} W ](\omega)
\\ &
- \int_{0}^{t} W(s) \Big(
\sqrt{1-|Z(s,\omega)|^{2}} - 1 \Big) e^{-i\omega s} ds e^{i \omega t}
\\ &
+ \int_{0}^{t} w(s)
\displaystyle
\frac{\Re [ \overline{Z(s,\omega)} \xi(s,\omega) ]}{\sqrt{1-|Z(s,\omega)|^{2}}}
e^{-i\omega s} ds e^{i \omega t},
\end{array}$$
and working as in (\ref{calcul-type-L2}), we get
\begin{multline*}\|\xi(t)\|_{L^{2}(\mathbb{R})}
\leqslant \sqrt{2\pi}
\| W \|_{L^{2}(0,T)} +
\| W \|_{L^{2}(0,T)} c \sqrt{T} \|Z\|_{C^{0}([0,T],L^{2}(\R))}
\\
+
\| w \|_{L^{2}(0,T)} c \sqrt{T} \| \xi \|_{C^{0}([0,T],L^{2}(\R))}.
\end{multline*}
Using (\ref{Z-borne-L-2}) and $c \sqrt{T} R \leqslant 1/2$, we get (\ref{diff-cont}).
\\

\emph{Third step: $F_T$ is a closed map because $6 c \sqrt{T} R \leqslant 1/2$.}
\\
Let $A$ be a closed subset of $B_R[L^{2}(0,T)]$.
Let $(Z_{n}(T,.)=F_T(w_n))_{n \in \mathbb{N}}$ be a sequence of $F_T(A)$ that
converges to $Z_{\infty}(.)$ in $L^{2} \cap C^{0}_{b}(\mathbb{R})$.
In order to prove that $Z_{\infty} \in F_T(A)$,
we prove that $(w_n)_{n \in \mathbb{N}}$ is a Cauchy
sequence in $L^{2}(0,T)$. For every $n \in \mathbb{N}$, we have
$$Z_{n}(T,\omega)= - \mathcal{F} [ \tau_{-T} w_n ](\omega)
- \int_{0}^{T} w_n(t) \Big( \sqrt{1-|Z_n(t,\omega)|^{2}} - 1 \Big)
e^{-i\omega t} dt e^{i \omega T},$$
so, for $n, p \in \mathbb{N}$, we have
$$\begin{array}{ll}
\mathcal{F} [ \tau_{-T} (w_n-w_p) ](\omega) =
&
(Z_p-Z_n)(T,\omega)
\\ &
- \int_0^T (w_n-w_p)(t) \Big( \sqrt{1-|Z_n(t,\omega)|^{2}} - 1  \Big)
e^{-i\omega t} dt e^{i \omega T}
\\ &
- \int_0^T w_p(t) \Big(
\sqrt{1-|Z_n(t,\omega)|^{2}} - \sqrt{1-|Z_p(t,\omega)|^{2}}
\Big) e^{-i\omega t} dt e^{i \omega T}.
\end{array}$$
Considering the $L^{2}(\mathbb{R})$-norm of each side,
using the Plancherel equality and working as in (\ref{calcul-type-L2}), we get
$$\begin{array}{ll}
\sqrt{2\pi} \| w_n-w_p \|_{L^{2}(0,T)} \leqslant
&
\|(Z_n-Z_p)(T)\|_{L^{2}(\mathbb{R})}
\\ &
+ \|w_n-w_p\|_{L^{2}(0,T)} c \sqrt{T} \|Z_n\|_{C^{0}([0,T],L^{2}(\R))}
\\ &
+ \|w_p\|_{L^{2}(0,T)} c \sqrt{T} \| Z_n - Z_p \|_{C^{0}([0,T],L^{2}(\R))}.
\end{array}$$
Using (\ref{Z-borne-L-2}) and (\ref{continuity-L-2}), we get
$$\sqrt{2\pi}\| w_n - w_p \|_{L^{2}(0,T)} \leqslant
\|(Z_n-Z_p)(T)\|_{L^{2}(\mathbb{R})}
+  6 c \sqrt{2\pi}\sqrt{T} R  \|w_n-w_p\|_{L^{2}(0,T)}.$$
which gives the conclusion because $6 c \sqrt{T} R = 1/2$.
\\

\emph{Fourth step: The manifold  $S:=F_T { B_R [ L^{2}(0,T) ] }$
does not coincide with its tangent space at zero.}
\\
We have
$$dF_{T}(0).W=- \mathcal{F}[ \tau_{-T} W ]$$
thus
$$T_{0}S=\mathcal{F}[ L^{2}((-T,0)) ].$$

Let us compute the third order development of $F_T$ around $0$,
$$w(t)= \epsilon W_1(t)   + \epsilon^2 W_2(t)   + \epsilon^3 W_3(t)   + ...$$
$$Z(t,\omega)= \epsilon Z_1(t,\omega) + \epsilon^2 Z_2(t,\omega)
+ \epsilon^3 Z_3(t,\omega) + ...$$
Since, as $\epsilon \rightarrow 0$,
$$\sqrt{1 - |\epsilon Z_1 + \epsilon^2 Z_2 + \epsilon^3 Z_3 |^{2}}
=
\sqrt{1 - \epsilon^{2} |Z_1|^{2} + o (\epsilon^{2}) }
=
1 - \frac{1}{2}\epsilon^{2} |Z_1|^{2} + o (\epsilon^{2}),$$
we have
$$Z_1(t,\omega)= - \mathcal{F}[ \tau_{-t} (W_1)_{|[0,t]} ](\omega),$$
$$Z_2(t,\omega) = - \mathcal{F}[ \tau_{-t}(W_2)_{|[0,t]} ](\omega),$$
$$Z_3(t,\omega) = - \mathcal{F}[ \tau_{-t}(W_3)_{|[0,t]} ](\omega)
- \frac{1}{2} \int_{0}^{t} W_1(\tau) |Z_1(\tau,\omega)|^{2}
e^{-i\omega \tau} d\tau e^{i \omega t}.$$
We want to prove the existence of $W_1 \in L^{2}(0,T)$ such that the map
$$\omega \mapsto \int_{0}^{T} W_1(\tau) |Z_1(\tau,\omega)|^{2}
e^{-i\omega \tau} d\tau$$
does not belong to $\mathcal{F}[L^{2}(0,T)]$.
Using the explicit expression of $Z_1$,
 the change of variable $\sigma_2 \rightarrow x=\tau+\sigma_1-\sigma_2$
and Fubini's theorem, we get
$$\begin{array}{rcl}
\int_{0}^{T} W_1(\tau) |Z_1(\tau,\omega)|^{2}  e^{-i\omega \tau} d\tau \phantom{\hspace{-2cm}}&&
\\ &=&
\int_{0}^{T} W_1(\tau)
\Big( \int_{0}^{\tau} W_1(\sigma_1) e^{-i\omega \sigma_1} d \sigma_1 \Big)
\Big( \int_{0}^{\tau} \overline{W_1(\sigma_2)} e^{i\omega \sigma_2} d \sigma_2 \Big)
e^{-i\omega \tau} d\tau
\\ &=&
\int_{\tau=0}^{T} \int_{\sigma_1=0}^{\tau} \int_{\sigma_2=0}^{\tau}
W_1(\tau)  W_1(\sigma_1)  \overline{W_1(\sigma_2)}
e^{-i \omega ( \tau+\sigma_1-\sigma_2)}
d\sigma_2 d\sigma_1 d\tau
\\ &=&
\int_{\tau=0}^{T} \int_{\sigma_1=0}^{\tau} \int_{x=\sigma_1}^{\sigma_1+\tau}
W_1(\tau)  W_1(\sigma_1)  \overline{W_1(\tau+\sigma_1-x)}  e^{-i\omega x}
dx d\sigma_1 d\tau
\\ &=&
\mathcal{F}[ \Phi_{W_1} ] (\omega)
\end{array}$$
where
$$\Phi_{W}(x):=\left\lbrace \begin{array}{l}
\int_{\tau=0}^{T} \int_{\sigma_1=\max\{ 0 , x-\tau \}}^{\min \{ \tau , x \}}
W(\tau)  W(\sigma_1)  \overline{W(\tau+\sigma_1-x)}
d\sigma_1 d\tau \text{ if } x \in (0,2T) \\
0 \text{ if } x \notin (0,2T).
\end{array} \right.$$
Computing the map $\Phi_1$ associated to $W = 1_{[0,T]}$,
we get $\Phi_1(3T/2)=T^2/16$ thus, for $\epsilon$ small enough,
$ F_T(\epsilon 1_{[0,T]})  \notin T_{0}S$. $\Box$

\subsection{Case $-\infty < \omega_* < \omega^* <+\infty$ }

The goal of this section is the proof of the following result.

\begin{thm} \label{analytic}

\textbf{(i)} Let $T>0$, $u,v \in L^2(0,T)$ and $M$ be the solution of
\begin{equation} \label{Bloch-sur-C}
\left\lbrace\begin{array}{l}
\frac{\partial M}{\partial t}(t,\omega)=
[\omega \Omega_z + u(t) \Omega_x + v(t) \Omega_y] M(t,\omega),\quad
(t,\omega) \in (0,T) \times \mathbb{C}, \\
M(0,\omega)=e_3.
\end{array}\right.
\end{equation}
Then $\omega \in \mathbb{C} \mapsto Z(T,\omega)$ is holomorphic.

\textbf{(ii)} Let $T>0$ and $R:=1/(4\sqrt{3T})$. There exists
$Z_f:(\omega_*,\omega^*)\rightarrow \mathbb{C}$ analytic such that,
for every $\epsilon^*>0$, there exists $\epsilon \in (0,\epsilon^*)$ such that,
for every $w \in B_R[L^2(0,T)]$, the solution of
(\ref{Bloch-ABC}) satisfies $Z(T) \neq \epsilon Z_f$.
\end{thm}

As a consequence, there are arbitrarily small analytic targets on $(\omega_*,\omega^*)$
that cannot be reached exactly in finite time,
with controls having a prescribed $L^2$-bound.

\textbf{Proof of Theorem \ref{analytic}:}
\textbf{(i).} Let $T>0$, $u,v \in L^2(0,T)$ and
$M$ be the solution of (\ref{Bloch-sur-C}).
We introduce the functions
$M_1,M_2:[0,T] \times \mathbb{R} \times \mathbb{R} \rightarrow \mathbb{R}^3$
defined by
$$\begin{array}{ll}
M_1(t,\omega_1,\omega_2):=\Re[ M(t,\omega_1 + i \omega_2)],
&
M_2(t,\omega_1,\omega_2):=\Re[ M(t,\omega_1 + i \omega_2)].
\end{array}$$
The function $\widetilde{M}(t,\omega_1,\omega_2):=
(M_1(t,\omega_1,\omega_2),M_2(t,\omega_1,\omega_2))^t$ solves an
equation of the form
\begin{equation} \label{Bloch-reel}
\frac{\partial \widetilde{M}}{\partial t}=f(t,\widetilde{M},\omega_1,\omega_2).
\end{equation}
The map $f$ is measurable, of class $C^1$ with respect to $(\widetilde{M},\omega_1,\omega_2) \in
\mathbb{R}^6 \times \mathbb{R} \times \mathbb{R}$ and satisfies
\begin{gather*}
|f(t,\widetilde{M},\omega_1,\omega_2)|\leqslant (|\omega_1|+|\omega_2|+|u|+|v|)|\widetilde{M}|
\\
|f_{\widetilde{M}}(t,\widetilde{M},\omega_1,\omega_2)|\leqslant |\omega_1|+|\omega_2|+|u|+|v|
\\
|f_{\omega_1}(t,\widetilde{M},\omega_1,\omega_2)|\leqslant |\widetilde{M}|,\,
|f_{\omega_2}(t,\widetilde{M},\omega_1,\omega_2)|\leqslant |\widetilde{M}|.
\end{gather*}
Thus
$\widetilde{M}$ has partial derivatives with respect to $\omega_1$ and $\omega_2$. (To check that, one can, for example, adapt the proof of
\cite[Theorem 4.1, chap. 4, page 100]{Hartman}.)
Let us prove that they satisfy the Cauchy-Riemann relations, in order to
get the holomorphy of $ \omega \in \mathbb{C} \mapsto M(T,\omega)$.
We introduce the notation
$$Y_{k,l}(t,\omega_1,\omega_2):=
\frac{\partial M_k}{\partial \omega_l}(t,\omega_1,\omega_2),
\text{ for } k,l \in \{1,2\}.$$
Differentiating the system (\ref{Bloch-reel}) with respect to $\omega_1$ and
$\omega_2$, we get
$$\left\lbrace\begin{array}{l}
\frac{\partial(Y_{11}-Y_{22})}{\partial t}=
A(t,\omega_1) (Y_{11}-Y_{22})
-
B(\omega_2)(Y_{12}+Y_{21}),
\\
\frac{\partial(Y_{12}+Y_{21})}{\partial t}=
A(t,\omega_1) (Y_{12}+Y_{21})
+
B(\omega_2) (Y_{11}-Y_{22}),
\\
(Y_{11}-Y_{22})(0,\omega_1,\omega_2)=0,
\\
(Y_{12}+Y_{21})(0,\omega_1,\omega_2)=0,
\end{array}\right.$$
where
$$A(t,\omega_1):=
\left( \begin{array}{ccc}
0        & - \omega_1 & v(t)  \\
\omega_1 &      0     & -u(t) \\
-v(t)    &     u(t)   & 0
\end{array} \right)
\text{ and }
B(\omega_2):=
\left( \begin{array}{ccc}
0        & - \omega_2 & 0  \\
\omega_2 &      0     & 0 \\
0        &     0      & 0
\end{array} \right).$$
The uniqueness of the solution of the Cauchy problem ensures that
$Y_{11}=Y_{22}$ and $Y_{12}=-Y_{21}$.
\\

\textbf{(ii)} Let $Z_f:\mathbb{R} \rightarrow \mathbb{C}$ be an analytic function
that does not belong to the tangent space of the image of $F_T$ at zero
(i.e. which is not the Fourier transform of an a function in $L^2((-T,0))$).
Then, for every $\epsilon^*>0$, there exists $\epsilon \in (0,\epsilon^*)$ such
that $\epsilon Z_f$ does not belong to the image of $F_T$.
Thanks to \textbf{(i)}, reaching $\epsilon Z_f$ on $(\omega_*,\omega^*)$
in time $T$ with controls $u$ and $v$ in $L^2((0,T),\R)$ is equivalent to reaching it on
$\mathbb{R}$. But $\epsilon Z_f$
does not belong to the image of $F_T$. Therefore $\epsilon Z_f$ cannot
be reached on $(\omega_*,\omega^*)$, in time $T$,
with  controls $u$ and $v$ in $B_R[L^2((0,T),\R)]$. $\Box$

\section{Ensemble controllability with unbounded controls}
\label{Li-Khaneja}

In this section, we take $-\infty < \omega_* < \omega^* < +\infty$.
The goal of this section is to complete the
very interesting arguments of ~\cite{Li-Khaneja:cdc06,Li-Khaneja:PRA06,Li-Khaneja:ieee09}
with functional analysis ideas, to prove the ensemble controllability of (\ref{Bloch-ABC}),
(i.e. the approximate controllability of (\ref{Bloch-ABC}) in $L^2(\omega_*,\omega^*)$)
with unbounded controls. Actually, we prove a stronger result.
\\

 First, let us introduce the definition of solutions of (\ref{Bloch-ABC}) with
Dirac controls.

\begin{Def}
\label{defDirac-Cauchy-Bloch}
Let $b \in [0,+\infty)$, $\beta, \gamma \in \mathbb{R}$ and
$M_0:(\omega_*,\omega^*) \rightarrow \mathbb{S}^2$. The solution of
(\ref{Bloch-ABC}) with $M(0)=M_0$, $u(t)=\beta \delta_b(t)$,
$v(t)=\gamma \delta_b(t)$ is
$$M(t,\omega)=\left\lbrace \begin{array}{l}
\exp( \omega \Omega_z t ) M_0(\omega) \text{ for } t \in [0,b),
\\
\exp( \omega \Omega_z (t-b) ) \exp( \beta \Omega_x + \gamma \Omega_y)
\exp( \omega \Omega_z b) M_0(\omega) \text{ for } t \in (b,+\infty),
\end{array}\right.$$
i.e.
$$M(b^+,\omega)=\exp( \beta \Omega_x + \gamma \Omega_y)M(b^-,\omega).$$
\end{Def}

Let
$$H^1((\omega_*,\omega^*),\mathbb{S}^2) :=
\{ M \in H^1((\omega_*,\omega^*),\mathbb{R}^3) ;
M(\omega) \in \mathbb{S}^2 , \forall \omega \in (\omega_*,\omega^*) \}.$$
Let also $U[t ; u , v ,M_0]$ denote the value at time $t$,
of the solution of (\ref{Bloch-ABC}), with initial condition $M_0$ at time $0$.
Thus, $U[t ; u , v ,M_0]$ is a function of $\omega \in (\omega_*,\omega^*)$. Definition
\ref{defDirac-Cauchy-Bloch} is motivated by the following result,
which is a consequence of explicit expressions and the
boundedness of $(\omega_*,\omega^*)$.

\begin{Prop} \label{cv:dirac}
 For every $\beta, \gamma \in \mathbb{R}$,  we have
$$\lim_{\epsilon \rightarrow 0}\Big\| U \Big[ b+\epsilon ; \frac{\beta}{\epsilon}1_{[b,b+\epsilon]},
\frac{\gamma}{\epsilon}1_{[b,b+\epsilon]},.\Big] -
U[b^+;\beta \delta_b,\gamma \delta_b,.]
\Big\|_{\mathcal{L}(H^{1}((\omega_*,\omega^*),\R^3), H^{1}((\omega_*,\omega^*),\R^3))}
=0.$$
\end{Prop}

Let us introduce the set $D$ of finite sums of Dirac masses on $[0,+\infty)$.
The goal of this section is to  prove the following result.

\begin{thm} \label{thm:LK}
Let $M_0 \in H^1((\omega_*,\omega^*),\mathbb{S}^2)$.
There exist $(t_n)_{n \in \mathbb{N}} \in [0,+\infty)^{\mathbb{N}}$,
and $(u_n)_{n \in \mathbb{N}}, (v_{n})_{n \in \mathbb{N}} \in D^{\mathbb{N}} $
such that
$$U[t_n^+;u_n,v_n,M_0] \rightarrow e_3
\text{ weakly in } H^1((\omega_*,\omega^*),\mathbb{R}^3).$$
\end{thm}

Thanks to Proposition \ref{cv:dirac}, one easily get, from Theorem \ref{thm:LK}
the following corollary.

\begin{Cor} \label{corollaire}
Let $M_0 \in H^1((\omega_*,\omega^*),\mathbb{S}^2)$.
There exist $(t_n)_{n \in \mathbb{N}} \in [0,+\infty)^{\mathbb{N}}$,
and $(u_n)_{n \in \mathbb{N}}, (v_{n})_{n \in \mathbb{N}} \in
L^{\infty}_{loc}([0,+\infty),\mathbb{R})^{\mathbb{N}}$
such that
$$U[t_n;u_n,v_n,M_0] \rightarrow e_3
\text{ weakly in } H^1((\omega_*,\omega^*),\mathbb{R}^3).$$
\end{Cor}

Thanks to the compactness of the injection
$H^{1}(\omega_*,\omega^*) \rightarrow L^2(\omega_*,\omega^*)$,
and the time reversibility of (\ref{Bloch-ABC}),
Theorem \ref{thm:LK} and Corollary \ref{corollaire}
give the ensemble controllability of (\ref{Bloch-ABC})
(i.e. its approximate controllability,
for the $L^2((\omega_*,\omega^*),\mathbb{R}^3)$-norm, in finite time).
These statements also give the approximate controllability of (\ref{Bloch-ABC}),
for the $L^{\infty}((\omega_*,\omega^*),\mathbb{R}^3)$-norm,
in finite time.
The proof of Theorem \ref{thm:LK} relies on the following Lemma,
that will be proved later on.

\begin{Lem} \label{Lem:dec}

\textbf{(1)} Let $M \in H^1((\omega_*,\omega^*),\mathbb{S}^2)$ be such that $M' \neq 0$.
There exist $T>0$, $u,v \in D$ such that
\begin{itemize}
\item one has
$$\| U[T^+;u,v,M]' \|_{L^2} < \| M' \|_{L^2},$$
\item for every sequence
$(M_n)_{n \in \mathbb{N}} \in H^1((\omega_*,\omega^*),\mathbb{S}^2)^{\mathbb{N}}$
satisfying
\begin{equation} \label{suitebornee}
\|M_n'\|_{L^2} \leqslant \| M' \|_{L^2}, \forall n \in \mathbb{N}
\end{equation}
and
\begin{equation} \label{suitefc}
M_{n} \rightarrow M \text{ weakly in } H^1((\omega_*,\omega^*),\mathbb{R}^3)
\end{equation}
there exists an extraction $\varphi$ such that
$$\| U[T^+;u,v,M_{\varphi(n)}]'  \|_{L^2} \leqslant \| M_{\varphi(n)}' \|_{L^2},
\forall n \in \mathbb{N}.$$
\end{itemize}

\textbf{(2)} Let $M \in \mathbb{S}^2$ be such that $M \neq e_3$.
There exists $\theta \in [0,2\pi)$ such that
for some $(u,v) \in \{
(\pi \delta_{1} + (\pi+\theta) \delta_{2},0) ,
( 0 , \pi \delta_{1} + (\pi+\theta) \delta_{2}) \}$,
$U[2^+;u,v, M]$ is constant over $(\omega_*,\omega^*)$ and
$$| U[2^+;u,v, M] - e_3 | < |M-e_3|.$$
\end{Lem}

In section \ref{subsec:Lem->Thm}, we prove
Theorem \ref{thm:LK} thanks to functional analysis and Lemma \ref{Lem:dec},
which is proved in section \ref{subsec:LKargument} and \ref{subsec:Proof of Lem:dec}.
In section \ref{subsec:LKargument}, we recall a preliminary result,
which is already presented in ~\cite{Li-Khaneja:cdc06,Li-Khaneja:ieee09}.
In section \ref{subsec:Proof of Lem:dec},
we deduce the proof of Lemma \ref{Lem:dec}.

\subsection{Proof of Theorem \ref{thm:LK} thanks to Lemma \ref{Lem:dec}}
\label{subsec:Lem->Thm}

In this section, we deduce Theorem \ref{thm:LK} from Lemma \ref{Lem:dec},
using similar arguments as in \cite{Nersesyan}.
\\

\textbf{Proof of Theorem \ref{thm:LK} thanks to Lemma \ref{Lem:dec}:}
Let $M_0 \in H^1((\omega_*,\omega^*), \mathbb{S}^2)$ be such that
$M_0 \neq e_3$ (otherwise $t_n \equiv 0$ gives the conclusion).
We introduce the set
$$\begin{array}{ll}
K :=
&
\{ \widetilde{M} \in H^1((\omega_*,\omega^*),\mathbb{S}^2) ;\,
\exists (t_n)_{n \in \mathbb{N}} \in  [0,\infty)^{\mathbb{N}},
\exists (u_n)_{n \in \mathbb{N}}, (v_n)_{n \in \mathbb{N}} \in  D^{\mathbb{N}}
\\ &
\text{ such that }
\| U[t_n^+;u_n,v_n,M_0]'\|_{L^2} \leqslant \| M_0' \|_{L^2}, \forall n \in \mathbb{N}
\\ & \text{ and }
U[t_n^+;u_n,v_n,M_0] \rightarrow \widetilde{M} \text{ weakly in }
H^1((\omega_*,\omega^*),\mathbb{R}^3) \}
\end{array}$$
and the quantity
$$m:= \inf \{ \|\widetilde{M}'\|_{L^2} ; \, \widetilde{M} \in K \}.$$
Notice that $K$ is not empty because it contains $M_0$ (take $t_n \equiv 0$).
\\

\emph{First step: Let us prove the existence of $e \in K$ such that
$\|e'\|_{L^2}=m$.}

Let $(M_n)_{n \in \mathbb{N}^*} \in K^{\mathbb{N}^*}$ be such that
$\| M_n' \|_{L^2} \rightarrow m$ when $n \rightarrow + \infty$.
Then $(M_n)_{n \in \mathbb{N}^*}$ is a bounded sequence in
$H^1((\omega_*,\omega^*),\mathbb{R}^3)$,
thus, there exists $e \in H^1((\omega_*,\omega^*),\mathbb{S}^2)$ such that (up to an extraction)
\begin{equation} \label{rel1.1}
M_{n} \rightarrow e \text{ weakly in } H^1
\text{ and strongly in } L^{2}.
\end{equation}
Then,
$$\|e'\|_{L^2} \leqslant
\liminf_{n \rightarrow + \infty} \| M_n' \|_{L^2} = m.$$
Let us prove that $e$ belongs to $K$, which gives the conclusion of the first step.
 For every $n \in \mathbb{N}^*$, $M_n \in K$ so
there exist $(t_{n}^p)_{p \in \mathbb{N}} \in [0,+\infty)^{\mathbb{N}}$,
$(u_{n}^{p})_{p \in \mathbb{N}}, (v_{n}^{p})_{p \in \mathbb{N}} \in D^{\mathbb{N}}$
such that
\begin{equation} \label{borneH1np}
\| U[t_n^{p+};u_n^p,v_n^p,M_0]'  \|_{L^2} \leqslant
\| M_0' \|_{L^2} , \forall p \in \mathbb{N}, \forall n \in \mathbb{N}^*
\end{equation}
and
$$U[t_n^{p+};u_n^p,v_n^p,M_0] \rightarrow M_n
\text{ weakly in } H^1
\text{ and strongly in } L^{2},
\text{ when } p \rightarrow + \infty,
\forall n \in \mathbb{N}^*.$$
 For every $n \in \mathbb{N}^*$, we choose $p=p(n) \in \mathbb{N}$ such that
\begin{equation} \label{rel1.2}
\| U[t_n^{p(n)+};u_n^{p(n)},v_n^{p(n)},M_0] - M_n \|_{L^{2}} \leqslant \frac{1}{n}.
\end{equation}
The sequence $( Y_n := U[t_n^{p(n)+};u_n^{p(n)},v_n^{p(n)},M_0] )_{n \in \mathbb{N}^*}$
is bounded in $H^1((\omega_*,\omega^*),\mathbb{R}^3)$ because of (\ref{borneH1np}).
Thus, there exists $e_{\sharp} \in H^1((\omega_*,\omega^*),\mathbb{S}^2)$
such that (up to an extraction)
$$ Y_n \rightarrow e_{\sharp}
\text{ weakly in } H^1
\text{ and strongly in } L^{2}.$$
The definition of $K$ ensures that $e_{\sharp}$ belongs to $K$.
Moreover, because of (\ref{rel1.1}) and (\ref{rel1.2}),
$Y_n \rightarrow e$
strongly in $L^2((\omega_*,\omega^*),\mathbb{R}^3)$, thus
(uniqueness of the strong $L^{2}$ limit), $e=e_{\sharp}$ and $e \in K$.
\\

\emph{Second step: Let us prove that $m=0$.}

Working by contradiction, we assume that $m>0$.
Then $e' \neq 0$ thus, we can apply Lemma \ref{Lem:dec} \textbf{(1)}.
There exist $T>0$, $u,v \in D$ such that
$$\| U[T^+;u,v,e]' \|_{L^2} < \|e'\|_{L^2}=m,$$
and an extraction $\varphi$ such that, with the notations of the first step,
\begin{equation} \label{U[T;u,v,Yn]bornee}
\| U[T^+;u,v,Y_{\varphi(n)}]' \|_{L^2} \leqslant
\|Y_{\varphi(n)}'\|_{L^2}, \forall  n \in \mathbb{N}.
\end{equation}
Let us prove that $U[T^+;u,v,e]$ belongs to $K$, which gives the contradiction.

Using (\ref{borneH1np}) and (\ref{U[T;u,v,Yn]bornee}), we have
$$\| U[T^+;u,v,Y_{\varphi(n)}]' \|_{L^2} \leqslant
\| M_0'  \|_{L^2}, \forall  n \in \mathbb{N}.$$
Thus, there exists $e_{*} \in H^{1}((\omega_*,\omega^*),\mathbb{S}^2)$ such that (up to an extraction)
$$U[T^+;u,v,Y_{\varphi(n)}] \rightarrow e_{*}
\text{ weakly in } H^1
\text{ and strongly in } L^{2}.$$
Then, $e_{*} \in K$, by definition of $K$.
But we have
$$\| U[T^+;u,v,Y_{\varphi(n)}] - U[T^+;u,v,e] \|_{L^2}
= \|Y_{\varphi(n)} - e \|_{L^2} \rightarrow 0 \text{ when } n \rightarrow + \infty,$$
thus $U[T^+,u,v,e]=e_{*}$ (uniqueness of the strong $L^{2}$ limit),
and $U[T^+,u,v,e] \in K$.  This ends the proof of the second step.
\\

\emph{Third step.} With a slight abuse of notations, let us still denotes by
$\mathbb{S}^2$ the set of constant functions from $(-\omega_*,\omega^*)$ with values into $\mathbb{S}^2$. Thanks to the first and second steps, the set $K \cap \mathbb{S}^2$ is not empty,
so we can consider
$$\widetilde{m}:=\inf \{ |\widetilde{M}-e_3| ; \widetilde{M} \in K \cap \mathbb{S}^2 \}.$$
Working exactly as in the first step, one can prove that
$K \cap \mathbb{S}^2$ is a closed subset of $\mathbb{S}^2$, thus $K \cap \mathbb{S}^2$ is compact
and there exists $\tilde{e} \in K \cap \mathbb{S}^2$ such that $|\tilde{e}-e_3| = \widetilde{m}$.
\\

\emph{Fourth step: Let us prove that $\widetilde{m}=0$, which gives the conclusion.}

Working by contradiction, we assume $\widetilde{m}>0$.
Then $\tilde{e} \neq e_3$ and we can apply Lemma \ref{Lem:dec} \textbf{(2)}.
There exists $\theta \in [0,2\pi)$ such that,
for some $(u,v) \in \{ (\pi \delta_{1} + (\pi+\theta) \delta_{2},0) ,
(0,\pi \delta_{1} + (\pi+\theta) \delta_{2}) \}$,
$U[2^+;u,v, \tilde{e}]$ is constant over $(\omega_*,\omega^*)$ and
$$| U[2^+;u,v, \tilde{e}] - e_3 | < |\tilde{e}-e_3|.$$
Let us prove that $U[2^+;u,v, \tilde{e}]$ belongs to $K$,
which gives the contradiction.

 First, we emphasize that explicit computations show that,
\begin{equation} \label{efface-drift}
\exp( \pi \Omega_{\xi} ) \exp( \omega \Omega_z ) \exp( \pi \Omega_{\xi})=
\exp( - \omega \Omega_z ), \forall \xi \in \{x,y\}.
\end{equation}
Thus, for some $\xi \in \{x,y\}$, we have
$$\begin{array}{ll}
U[2^+;u,v,.]
&
=\exp( (\theta + \pi) \Omega_{\xi} ) \exp ( \omega \Omega_z )
\exp(\pi \Omega_{\xi}) \exp( \omega \Omega_z )
\\ &
= \exp(\theta \Omega_{\xi})
\exp(\pi \Omega_{\xi}) \exp ( \omega \Omega_z ) \exp(\pi \Omega_{\xi})
\exp( \omega \Omega_z )
\\ &
=  \exp(\theta \Omega_{\xi}) \exp( - \omega \Omega_z ) \exp( \omega \Omega_z )
\\ &
= \exp( \theta \Omega_{\xi}).
\end{array}$$
Since $\tilde{e} \in K$, there exist
$(s_n)_{n \in \mathbb{N}} \in [0,+\infty)^{\mathbb{N}}$,
$(\mu_n)_{n \in \mathbb{N}}, (\nu_n)_{n \in \mathbb{N}} \in D^{\mathbb{N}}$
such that
$$\| U[s_n^+;\mu_n,\nu_n,M_0]'\|_{L^2} \leqslant \|M_0'\|_{L^2},
\forall n \in \mathbb{N},$$
$$U[s_n^+;\mu_n,\nu_n,M_0] \rightarrow \tilde{e}
\text{ weakly in } H^1.$$
Let $Z_n:=U[s_n^+;\mu_n,\nu_n,M_0]$. We have
$$U[2^+;u,v,Z_n]=\exp( \theta \Omega_{\xi} )Z_n.$$
Thus
$$\| U[2^+;u,v,Z_n]' \|_{L^2} = \| Z_n' \|_{L^2} \leqslant \|M_0'\|_{L^2},
\forall n \in \mathbb{N}$$
and $U[2^+;u,v,Z_n] \rightarrow \exp(\theta \Omega_{\xi})\tilde{e}=
U[2^+;u,v, \tilde{e}]$ weakly in $H^1$. Thus
$U[2^+;u,v, \tilde{e}] \in K$.$\Box$

\subsection{The argument of~\cite{Li-Khaneja:cdc06,Li-Khaneja:ieee09}}
\label{subsec:LKargument}

The goal of this section is to recall the proof of the following result,
which is already presented in ~\cite{Li-Khaneja:cdc06,Li-Khaneja:ieee09}.

\begin{Prop} \label{Prop:LK}
Let $P,Q \in \mathbb{R}[X]$. The flow of (\ref{Bloch-ABC}) can generate
$$I+\tau [P(\omega)\Omega_x+Q(\omega) \Omega_y] + o(\tau)
\text{ when } \tau \rightarrow 0,$$
with controls that are finite sums of Dirac masses.
More precisely, for every $\epsilon>0$,
there exists $\tau^*=\tau^*(P,Q,\epsilon)>0$ such that,
for every $\tau \in [0,\tau^*]$,
there exist $T>0$ and $u,v \in D$, such that
$$\Big\|  U[T^+;u,v,.] -  \Big( I+\tau [P(\omega)\Omega_x+Q(\omega) \Omega_y] \Big)
\Big\|_{\mathcal{L}(H^{1}((\omega_*,\omega^*),\R^3), H^{1}((\omega_*,\omega^*),\R^3))}
\leqslant
\epsilon \tau.$$
\end{Prop}

\begin{rk}
Let us explain why Proposition \ref{Prop:LK} may not be sufficient to prove the
ensemble controllability (i.e. the approximate $L^2(\omega_*,\omega^*)$-controllability)
of (\ref{Bloch-ABC}) with Dirac controls.

 First, let us remark that, for every point $M=(x^{(1)},x^{(2)},x^{(3)}) \in \mathbb{S}^2$,
such that $x^{(3)} \neq 0$, then $(\Omega_x M, \Omega_y M)$ is a basis of
$T_{\mathbb{S}^2}M$ (the tangent space of $\mathbb{S}^2$ at $M$).

Let $\epsilon>0$ and $M_0 = (x_0,y_0,z_0) \in L^2((\omega_*,\omega^*), \mathbb{S}^2)$ be such that,
$z(\omega) \neq 0$ for almost every $\omega \in (\omega_*,\omega^*)$.
 Following a classical strategy, we  consider an homotopy
$$\begin{array}{cccccc}
H: & [0,1] & \times & (\omega_*,\omega^*) & \rightarrow & \mathbb{S}^2 \\
   &  (s   &    ,   &   \omega)           & \mapsto     & H(s,\omega)
\end{array}$$
such that $H \in C^1([0,1],L^{2}((\omega_*,\omega^*),\mathbb{S}^2))$,
$H(0,\omega)=M_0(\omega)$, $H(1,\omega)=e_3$,
and we try to reach $e_3$ from $M_0$ by following the path given by $H$.
Since $z \neq 0$ a.e. on $(\omega_*,\omega^*)$, there exist
$f,g \in L^2((\omega_*,\omega^*),\mathbb{R})$ such that
$$\frac{\partial H}{\partial s}(0,\omega)=
f(\omega) \Omega_x M_0(\omega) + g(\omega) \Omega_y M_0(\omega).$$
Thanks to the Weierstrass theorem, there exist $P,Q \in \mathbb{R}[X]$
such that $\|f-P\|_{L^2} < \epsilon$ and $\|g-Q\|_{L^2} < \epsilon$.
Applying Proposition \ref{Prop:LK}, one may follow (approximately)
the direction given by $\frac{\partial H}{\partial s}(0,\omega)$,
with a small amplitude $\tau^*$, that depends on this direction.

If one wants to be sure to reach $e_3$ in finite time, by
iteration of this process, one would need, at least, the independence
of the amplitude $\tau^*$ with respect to the direction
(otherwise, one may stop in the middle of the path).
However the maximum amplitude $\tau^*$ given by Proposition \ref{Prop:LK}
depends on the direction, through the choice of the polynomials $P,Q$.
\end{rk}

\textbf{Proof of Proposition \ref{Prop:LK}:}
In this proof $\tau$ is a positive real number. By (\ref{efface-drift}),
\begin{equation} \label{cancel-drift-term}
\exp( \pi \Omega_x ) \exp( \omega \Omega_z \tau ) \exp( - \pi \Omega_x)=
\exp( - \omega \Omega_z \tau ),
\end{equation}
and this evolution is generated, in time $\tau$, by the controls
$u(t)=- \pi \delta_0(t) + \pi \delta_{\tau}(t)$, $v \equiv 0$.

The controls
$$u(t)=\sqrt{\tau} \delta_0(t) - (\pi + \sqrt{\tau}) \delta_{\sqrt{\tau}}(t)
+ \pi \delta_{2 \sqrt{\tau}}(t)$$
$$\text{(resp. }
u(t) = - \pi \delta_0(t) + (\pi-\sqrt{\tau}) \delta_{\sqrt{\tau}}(t)
+ \sqrt{\tau} \delta_{2 \sqrt{\tau}}(t) \text{)},$$
$v \equiv 0$, generate in time $2 \sqrt{\tau}$ the evolution
$$\begin{array}{ll}
U_1(\tau)
& :=
\exp(\pi \Omega_x) \exp(\omega \Omega_z \sqrt{\tau}) \exp(-(\pi+\sqrt{\tau})\Omega_x)
\exp(\omega \Omega_z \sqrt{\tau}) \exp(\Omega_x \sqrt{\tau})
\\ & =
\exp( - \omega \Omega_z \sqrt{\tau} ) \exp(- \Omega_x \sqrt{\tau} )
\exp(\omega \Omega_z \sqrt{\tau}) \exp(\Omega_x \sqrt{\tau})
\\ & =
I + \tau \omega [\Omega_z,\Omega_x] + o(\tau)
\\ & =
I + \tau \omega \Omega_y  + o(\tau)
\end{array}$$
(resp.
$$\begin{array}{ll}
U_1(-\tau)
& :=
\exp(\Omega_x \sqrt{\tau}) \exp(\omega \Omega_z \sqrt{\tau})
\exp((\pi-\sqrt{\tau})\Omega_x) \exp(\omega \Omega_z \sqrt{\tau})
\exp(-\pi \Omega_x)
\\ & =
\exp(\Omega_x \sqrt{\tau})\exp(\omega \Omega_z \sqrt{\tau})
\exp(- \Omega_x \sqrt{\tau} ) \exp( - \omega \Omega_z \sqrt{\tau} )
\\ & =
I - \tau \omega [\Omega_z,\Omega_x] + o(\tau)
\\ & =
I - \tau \omega \Omega_y  + o(\tau))
\end{array}$$
where we have used (\ref{cancel-drift-term}) to pass from the first to the second line. Here and in the following,  $o(\tau)$ denote quantities which tend to $0$ in the $\mathcal{L}(H^{1}((\omega_*,\omega^*),\R^3), H^{1}((\omega_*,\omega^*),\R^3))$-norm as $\tau \rightarrow 0^+$.
In the same way, there exist controls, that are sums of Dirac masses,
that generate in time $6 \sqrt{\tau}$ the evolutions
$$\begin{array}{ll}
U_2(\tau) &  :=
\exp( -\omega \Omega_z \sqrt{\tau} ) U_1(-\tau)
\exp(\omega \Omega_z \sqrt{\tau} ) U_1(\tau)
\\ & =
I + \tau^{3/2} \omega^2 [\Omega_y,\Omega_z] + o( \tau^{3/2} )
\\ & =
I - \tau^{3/2} \omega^2 \Omega_x + o( \tau^{3/2} ),
\end{array}$$
$$\begin{array}{ll}
U_2(-\tau) &  :=
U_1(\tau) \exp(\omega \Omega_z \sqrt{\tau} ) U_1(-\tau) \exp( -\omega \Omega_z \sqrt{\tau} )
\\ & =
I - \tau^{3/2} \omega^2 [\Omega_y,\Omega_z] + o( \tau^{3/2} )
\\ & =
I + \tau^{3/2} \omega^2 \Omega_x + o( \tau^{3/2} ),
\end{array}$$
and in time $14 \sqrt{\tau}$ the evolutions
$$\begin{array}{ll}
U_3(\tau) & :=
\exp( -\omega \Omega_z \sqrt{\tau} ) U_2(-\tau)
\exp(\omega \Omega_z \sqrt{\tau} ) U_2(\tau)
\\ & =
I - \tau^{2} \omega^3 [\Omega_z,\Omega_x] + o( \tau^{2} )
\\ & =
I - \tau^{2} \omega^3 \Omega_y + o( \tau^{2} ),
\end{array}$$
$$\begin{array}{ll}
U_3(-\tau) &  :=
U_2(\tau) \exp(\omega \Omega_z \sqrt{\tau} ) U_2(-\tau) \exp( -\omega \Omega_z \sqrt{\tau} )
\\ & =
I + \tau^{2} \omega^3 [\Omega_z,\Omega_x] + o( \tau^{2} )
\\ & =
I + \tau^{2} \omega^3 \Omega_y + o( \tau^{2} ),
\end{array}$$
Thus, one can generate
$I \pm  \tau \omega^2 \Omega_x + o(\tau)$ in time $6 \tau^{1/3}$,
and $I \pm \tau \omega^3 \Omega_y  + o(\tau)$ in time $14 \tau^{1/4}$.
Iterating this process, for every $n \in \mathbb{N}$, one can generate
$I \pm \tau \omega^{2n}     \Omega_x   + o(\tau)$ and
$I \pm \tau \omega^{2n+1} \Omega_y  + o(\tau)$ in a time
$T_n$ that behaves like $4^n \tau^{\frac{1}{2n}}$.
The same argument with $\Omega_x$ replaced by $\Omega_y$ in $U_{j}(\tau)$, $j \geqslant 1$,
shows that for every $n \in \mathbb{N}$, one can generate
$I \pm  \tau \omega^{2n+1} \Omega_x + o(\tau)$ and
$I \pm  \tau \omega^{2n} \omega_y  + o(\tau)$ in a time
$T_n$ that behaves like $4^n \tau^{\frac{1}{2n}}$.
Thus, for every $P,Q \in \mathbb{R}[X]$, one can generate
$$I + \tau [P(\omega) \Omega_x + Q(\omega) \Omega_y] + o(\tau)$$
in finite time, by composing the previous evolutions. $\Box$

\subsection{Proof of Lemma \ref{Lem:dec}}
\label{subsec:Proof of Lem:dec}

The goal of this subsection is the proof of Lemma \ref{Lem:dec},
thanks to the previous subsection.
\\

\textbf{Proof of Lemma \ref{Lem:dec}:}

\textbf{Proof of (1) of Lemma \ref{Lem:dec}.} It is sufficient
to prove this statement under the additional assumption
\begin{equation} \label{add_assump}
z \neq 0.
\end{equation}
Indeed, let us assume that it is proved
when (\ref{add_assump}) holds. Let $M=(x,y,z) \in H^1((\omega_*,\omega^*),\mathbb{S}^2)$
be such that $M' \neq 0$ and $z \equiv 0$. Then
$x^2 + y^2 \equiv 1$ thus $x \neq 0$ or $y \neq 0$. Let us assume, for example, that
$y \neq 0$. Thanks to (\ref{efface-drift}), we have
$$\begin{array}{ll}
U[2; \frac{3\pi}{2} \delta_0 + \pi \delta_1,0,.]
& =
\exp( \omega \Omega_z ) \exp (\pi \Omega_x) \exp(\omega \Omega_z )
\exp(\frac{3\pi}{2} \Omega_x)
\\ & =
\exp( \omega \Omega_z ) \exp (\pi \Omega_x) \exp(\omega \Omega_z )
\exp( \pi  \Omega_x) \exp (\frac{\pi}{2} \Omega_x)
\\ & =
\exp( \omega \Omega_z ) \exp( -\omega \Omega_z ) \exp (\frac{\pi}{2} \Omega_x)
\\ & =
\exp (\frac{\pi}{2} \Omega_x).
\end{array}$$
Thus the function
$$U[2; \frac{3\pi}{2} \delta_0 + \pi \delta_1,0,M]=
\left(\begin{array}{c}
x \\ 0 \\ y
\end{array}\right)$$
has a non vanishing third component and the $L^2$ norm of its derivative is
the same one as $M$. Applying Lemma \ref{Lem:dec} \textbf{(1)} to
$U[2; \frac{3\pi}{2} \delta_0 + \pi \delta_1,0,M]$, we get the conclusion
of Lemma \ref{Lem:dec} \textbf{(1)} for $M$.
\\

Let $M=(x,y,z) \in H^1((\omega_*,\omega^*),\mathbb{S}^2)$
be such that $M' \neq 0$ and $z \neq 0$.
\\

\emph{First step: Let us prove the existence of $P,Q \in \mathbb{R}[X]$,
$\alpha >0$ and $\tau^*_0>0$ such that, for every $\tau \in (0,\tau^*_0)$,
\begin{itemize}
\item one has
\begin{equation} \label{Lem:dec-step1}
\Big\| \frac{d}{d\omega} \Big[
\Big( I + \tau [P(\omega)\Omega_x+Q(\omega) \Omega_y] \Big) M \Big]
\Big\|_{L^2}^2    \leqslant \| M' \|_{L^2}^2 - \tau \alpha,
\end{equation}
\item for every sequence $(M_n)_{n \in \mathbb{N}}
\in H^1((\omega_*,\omega^*),\mathbb{S}^2)^{\mathbb{N}}$
satisfying (\ref{suitebornee}) and (\ref{suitefc}),
there exists an extraction $\varphi$ such that
\begin{equation} \label{Lem:dec-step1-suiteext}
\Big\| \frac{d}{d\omega} \Big[
\Big( I + \tau [P(\omega)\Omega_x+Q(\omega) \Omega_y] \Big) M_{\varphi(n)} \Big]
\Big\|_{L^2}^2    \leqslant
\| M_{\varphi(n)}' \|_{L^2}^2 - \tau \alpha,
\forall n \in \mathbb{N}.
\end{equation}
\end{itemize}}

Developing the square, we get, for $\tau \geqslant 0$ and $P,Q \in \mathbb{R}[X]$,
$$\Big\| \frac{d}{d\omega} \Big[
\Big( I + \tau [P(\omega)\Omega_x+Q(\omega) \Omega_y] \Big) M  \Big]
\Big\|_{L^2}^2
=
\| M' \|_{L^2}^2 + 2 \tau A(P,Q) + \tau^2 B(P,Q),$$
$$\Big\| \frac{d}{d\omega} \Big[
\Big( I + \tau [P(\omega)\Omega_x+Q(\omega) \Omega_y] \Big) M_n   \Big]
\Big\|_{L^2}^2 = \| M_n' \|_{L^2}^2 + 2 \tau A_n(P,Q) + \tau^2 B_n(P,Q),$$
where $A(P,Q), A_n(P,Q), B(P,Q), B_n(P,Q)$ are real constants.
Straightforward computations give
$$\begin{array}{ll}
A (P,Q)
& =
\int_{\omega_*}^{\omega^*}
\langle \frac{d}{d\omega} \Big[ [P(\omega)\Omega_x + Q(\omega) \Omega_y] M(\omega) \Big],
M'(\omega) \rangle d\omega
\\
& =
\int_{\omega_*}^{\omega^*}
\Big(  P' [ -z y' + y z' ] + Q' [  z x' - x z'  ] \Big) d\omega.
\end{array}$$
We look for $P,Q \in \mathbb{R}[X]$ such that $A(P,Q)<0$.
Since $A$ is a linear form in $(P,Q)$ it is sufficient to prove that $A \neq 0$.
Working by contradiction, we assume $A=0$. Thanks to the density of polynomials
in $L^2(\omega_*,\omega^*)$, we have
\begin{gather}
\label{eqxyz}
zy'-yz' = 0,
\,
zx'-xz' = 0.
\end{gather}
Let $I$ be a nonempty connected component of $\{\omega \in (\omega_*, \omega^*); \, z(\omega)\not =0\}$. Since $z \neq 0$, such a $I$ exists. By (\ref{eqxyz}), there exist $a, b \in \mathbb{R}$ such that
$x(\omega)=az(\omega)$ and $y(\omega)=bz(\omega)$, $\forall \omega \in I$.
Since $M$ takes values in $\mathbb{S}^2$, we have
$$1 = x(\omega)^2 + y(\omega)^2 + z(\omega)^2 = (a^2+b^2+1) z(\omega)^2.$$
This shows that $I=(\omega_*, \omega^*)$ and that $M$ is constant over $(\omega_*,\omega^*)$, which is in contradiction
with the assumption $M' \neq 0$.
Therefore, there exist $P,Q \in \mathbb{R}[X]$ such that $A(P,Q)<0$.
\\

 For every $n \in \mathbb{N}$, we have
$$A_n(P,Q)=\int_{\omega_*}^{\omega^*}
P'(-z_n y_n' + y_n z_n') + Q'(z_n x_n' - x_n z_n').$$
Thanks to  (\ref{suitefc}), there exists an extraction $\varphi$ such that
$$M_{\varphi(n)} \rightarrow M
\text{ weakly in } H^1 \text{ and strongly in } L^2.$$
Then $A_{\varphi(n)}(P,Q) \rightarrow A(P,Q)$ when $n \rightarrow + \infty$.
Thus, we can assume that
\begin{equation} \label{}
A_{\varphi(n)}(P,Q) < \frac{3}{4} A(P,Q), \forall n \in \mathbb{N}
\end{equation}
(otherwise take another extraction). We have
$$\begin{array}{ll}
\sqrt{B(P,Q)} & :=
\Big\| \frac{d}{d\omega} [P(\omega) \Omega_x + Q(\omega) \Omega_y]M \Big\|_{L^2}
\\ & \leqslant
\| P' \|_{L^2} + \| Q' \|_{L^2} + [ \|P\|_{L^{\infty}} + \|Q\|_{L^{\infty}} ] \|M'\|_{L^2},
\end{array}$$
$$\begin{array}{ll}
\sqrt{B_n(P,Q)} & :=
\Big\| \frac{d}{d\omega} [P(\omega) \Omega_x + Q(\omega) \Omega_y]M_n \Big\|_{L^2}
\\ & \leqslant
\| P' \|_{L^2} + \| Q' \|_{L^2} + [ \|P\|_{L^{\infty}} + \|Q\|_{L^{\infty}} ] \|M_n'\|_{L^2}
\\ & \leqslant
\| P' \|_{L^2} + \| Q' \|_{L^2} + [ \|P\|_{L^{\infty}} + \|Q\|_{L^{\infty}} ] \|M'\|_{L^2}.
\end{array}$$
Let $\tau^*_0=\tau^*_0(M)>0$ be such that
$$\tau^*_0 \Big(
\| P' \|_{L^2} + \| Q' \|_{L^2} + [ \|P\|_{L^{\infty}} + \|Q\|_{L^{\infty}} ] \|M'\|_{L^2}
\Big)^2 < \frac{|A(P,Q)|}{2}.$$
Then, for every $\tau \in [0,\tau^*_0]$, we have
(\ref{Lem:dec-step1}) and (\ref{Lem:dec-step1-suiteext}) with $\alpha:=-A(P,Q)$.
\\

\emph{Second step: Conclusion.}

Let $P, Q$ be as in the first step.
Let $\epsilon_1>0$ be such that
\begin{equation} \label{def:epsilon1}
\epsilon_1 \|M\|_{H^1} < \frac{\alpha}{2 \|M'\|_{L^2}}.
\end{equation}
Let $\tau^*=\tau^*(P,Q,\epsilon_1)$ be as in Proposition \ref{Prop:LK}
and $\tau^*_1:=\min\{ \tau^* , \tau^*_0 \}$.
Thanks to Proposition \ref{Prop:LK}, there exist $T>0$,
$u,v \in L^{\infty}_{loc}([0,+\infty),\mathbb{R})$ such that
\begin{equation} \label{applPropLKreg}
\Big\| U[T^+;u,v,.] -
\Big( I + \tau_1^* [P(\omega)\Omega_x+Q(\omega) \Omega_y] \Big)
\Big\|_{\mathcal{L}(H^1,H^1)}
\leqslant \epsilon_1 \tau^*_1.
\end{equation}
Then, using (\ref{applPropLKreg}), (\ref{Lem:dec-step1}) and (\ref{def:epsilon1}),
we get
$$\begin{array}{lll}
\|U[T^+;u,v,M]'\|_{L^2}
& \leqslant &
\Big\| \frac{d}{d\omega} \Big[ U[T;u,v,M] -
\Big( I + \tau [P(\omega)\Omega_x+Q(\omega) \Omega_y] \Big)M \Big] \Big\|_{L^2}
\\ & &
+ \Big\| \frac{d}{d\omega} \Big[
\Big( I + \tau [P(\omega)\Omega_x+Q(\omega) \Omega_y] \Big) M \Big]
 \Big\|_{L^2}
\\ & \leqslant &
\epsilon_1 \tau^*_1 \|M\|_{H^1}
+ \Big( \| M' \|_{L^2}^2 - \alpha \tau^*_1 \Big)^{1/2}
\\ & \leqslant &
\epsilon_1 \tau^*_1  \|M\|_{H^1}
+ \| M'\|_{L^2} - \frac{\alpha \tau^*_1}{2 \| M'\|_{L^2} }
\\ & < &
\| M'\|_{L^2}.
\end{array}$$
Similarly, we have
$$\begin{array}{lll}
\|U[T^+;u,v,M_{\varphi(n)}]'\|_{L^2}
& \leqslant
\epsilon_1 \tau^*_1 \|M_{\varphi(n)}\|_{H^1}
+ \| M_{\varphi(n)}'\|_{L^2} - \frac{\alpha \tau^*_1}{2 \| M_{\varphi(n)}'\|_{L^2} }
\\ & \leqslant
\epsilon_1 \tau^*_1  \|M\|_{H^1}
+ \| M_{\varphi(n)}'\|_{L^2} - \frac{\alpha \tau^*_1}{2 \| M' \|_{L^2} }
\\ & <
\| M_{\varphi(n)}'\|_{L^2}.
\end{array}$$
This ends the proof of the first statement of Lemma \ref{Lem:dec}.
\\

\textbf{Proof of (2) of Lemma \ref{Lem:dec}.} Let $M=(x,y,z) \in \mathbb{S}^2$ be such that $M \neq e_3$.
Then $x \neq 0$ or $y \neq 0$. We assume, for example, that $y \neq 0$.
Thanks to (\ref{efface-drift}), we have
$$\begin{array}{ll}
U[2^+;\pi \delta_1+(\pi+\theta)\delta_2,0,.]
& =
\exp((\pi+\theta)\Omega_x) \exp(\omega \Omega_z)
\exp(\pi \Omega_x) \exp( \omega \Omega_z)
\\ & =
\exp(\theta \Omega_x)
\exp(\pi \Omega_x) \exp(\omega \Omega_z) \exp(\pi \Omega_x)
\exp( \omega \Omega_z)
\\ & =
\exp(\theta \Omega_x) \exp( - \omega \Omega_z) \exp( \omega \Omega_z)
\\ & =
\exp(\theta \Omega_x).
\end{array}$$
Thus,
$$U[2^+;\pi \delta_1+(\pi+\theta)\delta_2,0,M]=
\left( \begin{array}{ccc}
1 & 0 & 0 \\
0 & \cos(\theta) & - \sin(\theta) \\
0 & \sin(\theta) & \cos(\theta)
\end{array} \right)M.$$
We get the conclusion with $\theta \in [0,2\pi)$ such that
$$U[2^+;\pi \delta_1+(\pi+\theta)\delta_2,0,M]=
\left( \begin{array}{c}
x \\ 0 \\ \sqrt{y^2 + z^2}
\end{array} \right). \Box$$

\section{Explicit controls for the asymptotic exact controllability to $e_3$}
\label{Section:asymptotic}

In this section, $\omega_*=0$, $\omega^*=\pi$.
We propose explicit controls realizing the asymptotic exact controllability to $-e_3$,
locally around $-e_3$.
\\

 First, let us introduce some notations.
 For a function $f:(-\pi,\pi) \rightarrow \mathbb{C}$,
we denote by $c_n(f)$ its Fourier coefficients and by $N(f)$ their $l^1$-norm:
$$\begin{array}{ll}
c_{n}(f):=\displaystyle \frac{1}{2\pi} \int_{-\pi}^{\pi} f(\omega) e^{-in\omega} d\omega,
&
N(f):=\sum\limits_{n \in \mathbb{Z}} |c_n(f)|.
\end{array}$$
 For a function $f:(0,\pi) \rightarrow \mathbb{C}$, we define
$$\begin{array}{ll}
c_{n}(f):=c_n(\tilde{f}), \forall n \in \mathbb{Z},
&
N(f):=N(\tilde{f}),
\end{array}$$
where $\tilde{f}:(-\pi,\pi) \rightarrow \mathbb{C}$,
$\tilde{f}(\omega):=f(|\omega|)$.
 For a vector valued map
$M=(x,y,z):(0,\pi) \rightarrow \mathbb{R}^3$, we define
$N(M):=N(x)+N(y)+N(z)$. Then, we have the following results.

\begin{Lem} \label{Prel}
 For every $f,g:[0,\pi] \rightarrow \mathbb{C}$ such that
$N(f), N(g)<\infty$, we have
\begin{equation} \label{N(fg)}
N(fg) \leqslant N(f) N(g).
\end{equation}
 For every $(x,y)\in L^1((0,\pi),\R^2)$, we have, with $Z:=x+iy$,
\begin{equation} \label{NxyZ}
\frac{1}{2} \Big( N(x)+N(y) \Big) \leqslant
N(Z) \leqslant N(x)+N(y).
\end{equation}
 For every $M:[0,\pi] \rightarrow \mathbb{S}^2$ such that $N(M)<+\infty$ and
$z(\omega) > 0, \forall \omega \in [-\pi,\pi]$, we have
\begin{equation} \label{N(z-1)}
N(z-1) \leqslant 2 N(Z)^2.
\end{equation}
If, moreover, $N(Z) \leqslant 1/4$, then
\begin{equation} \label{N(Z)measures}
\frac{1}{2} N(Z) \leqslant N(M-e_3) \leqslant 3N(Z).
\end{equation}
\end{Lem}

As a consequence, for a map
$M:[0,\pi] \rightarrow \mathbb{S}^2$ such that $N(Z) \leqslant 1/4$ and $z>0$,
the quantity $N(Z)$ measures the $N$-distance from $M$ to $e_3$.
\\

\textbf{Proof of Lemma \ref{Prel}:}
We have
$$N(fg) =
\sum\limits_{n \in \mathbb{Z}} | \sum\limits_{p \in \mathbb{Z}}
c_{n-p}(f) c_{p}(g) |
 \leqslant
\sum\limits_{p \in \mathbb{Z}} |c_{p}(g) |
\sum\limits_{n \in \mathbb{Z}} |c_{n-p}(f)|
= N(f) N(g).$$
The inequality (\ref{NxyZ}) is a consequence of the triangular inequality because
$N(Z)=N(x+iy)$ and $N(x)+N(y)=N((Z+\overline{Z})/2)+N((Z-\overline{Z})/2i)$.

Let $M:[0,\pi] \rightarrow \mathbb{S}^2$ be such that $N(M)<+\infty$ and $z > 0$.
We have
$$N(z-1)
= \frac{1}{2\pi} \int_{-\pi}^{\pi} 1 - \sqrt{1-|\tilde{Z}(\omega)|^2} d\omega +
\sum\limits_{n \in \mathbb{Z}-\{0\}}
\Big| \frac{1}{2\pi} \int_{-\pi}^{\pi} \sqrt{1-|\tilde{Z}(\omega)|^2} e^{-in\omega} d\omega \Big|.$$
Using $1 - \sqrt{1-x} \leqslant x, \forall x \in (0,1)$,
$\sqrt{1-x} = 1 + \sum_{p=1}^{\infty} \alpha_p x^p$,
that converges uniformly with respect to $x \in [0,\delta]$ when $\delta<1$
and where $\alpha_p<0$ for every $p \in \mathbb{N}^*$ and (\ref{N(fg)}), we get
$$\begin{array}{ll}
N(z-1)
& \leqslant
\|\tilde{Z}\|_{L^{\infty}(0,2\pi)}^2 -
\sum\limits_{p=1}^{\infty} \alpha_p \sum\limits_{n \in \mathbb{Z}-\{0\}}
\Big| \frac{1}{2\pi} \int_{-\pi}^{\pi} |\tilde{Z}(\omega)|^{2p} e^{-in\omega} d\omega \Big|
\\ & \leqslant
N(Z)^2 - \sum\limits_{p=1}^{\infty} \alpha_p N(|Z|^{2p})
\\ & \leqslant
N(Z)^2 - \sum\limits_{p=1}^{\infty} \alpha_p N(Z)^{2p}
\\ & \leqslant
N(Z)^2 + 1 - \sqrt{1-N(Z)^2}
\\ & \leqslant 2 N(Z)^2.
\end{array}$$
 Formula (\ref{N(Z)measures}) is a direct consequence of the previous inequalities.$\Box$
\\

The goal of this section is the proof of the following theorem.

\begin{thm} \label{Main-thm-(0,pi)}
There exists $\delta>0$ such that,
for every $M_0 : [0,\pi] \rightarrow \mathbb{S}^2$ with
$N[ Z_0 ] < \delta$ and $z_0 < -1/2$,
there exists $\epsilon=\epsilon(M_0)>0$ such that,
the solution of (\ref{Bloch-ABC}) with $M(0)=M_0$,
$$u(t):= \frac{\pi}{\epsilon} 1_{[k,k+\epsilon]}(t)
- \sum\limits_{p=1}^{2k-1} \Im \Big(  c_{-k+p} (Z_0)  \Big)
\frac{1}{\epsilon} 1_{[k+p,k+p+\epsilon]}(t)
+ \frac{\pi}{\epsilon} 1_{[3k,3k+\epsilon]}(t),$$
$$v(t):= - \sum\limits_{p=1}^{2k-1}
\Re \Big( c_{-k+p} (Z_0) \Big) \frac{1}{\epsilon} 1_{[k+p,k+p+\epsilon]}(t),$$
where $k=k(M_0) \in \mathbb{N}$ is such that
\begin{equation} \label{cdt-k}
\sum\limits_{|n|>k} |c_n(Z_0)| < \frac{N(Z_0)}{4},
\end{equation}
satisfies
$$N[ Z(3k+\epsilon) ] < \frac{N[Z_0]}{2},$$
$$z(3k+\epsilon) < -1/2.$$
\end{thm}

By iterating this process, we find
an increasing sequence $(t_n)_{n \in \mathbb{N}} \in [0,+\infty)^{\mathbb{N}}$
and two controls $u,v \in L^{\infty}_{loc}([0,+\infty),\mathbb{R})$ such that
$$N[Z(t_n)] < \frac{1}{2^n} N[Z_0].$$
Thus, $\|M(t_n)+e_3\|_{L^{\infty}} \rightarrow 0$ when $n \rightarrow + \infty$.
These explicit controls provide the exact asymptotic controllability to $e_3$.

In section \ref{subsec:Heuristic}, we present the heuristic
of the proof of Theorem \ref{Main-thm-(0,pi)}.
which is detailed in section \ref{subsec:proof of thm}.

\subsection{Heuristic}
\label{subsec:Heuristic}

Let us sketch the proof of Theorem \ref{Main-thm-(0,pi)}.
It is inspired by the return method, introduced in \cite{92mcss,96jmpa} and already used for the control of quantum systems in \cite{05Beauchard,2006-jfa} (for other applications see the book \cite{coron:book}).
It consists here in going close to $+e_3$ in order to delete the
main Fourier coefficients of the initial condition, and then
to move back to $-e_3$.
\par

Notice that, when $z>0$, the system (\ref{Bloch-ABC}) implies that
\begin{equation} \label{eq:Z-sym}
\dot{Z}(t,\omega)=i \omega Z(t,\omega) - w(t) \sqrt{1-|Z(t,\omega)|^2},
(t,\omega) \in (0,+\infty) \times (0,\pi),
\end{equation}
\begin{equation} \label{eq:z}
\dot{z}(t,\omega) = - \Re[ w(t) \overline{Z}(t,\omega)],
(t,\omega) \in (0,+\infty) \times (0,\pi).
\end{equation}

We have
\begin{equation} \label{dec:Z0}
Z_0(\omega)=\sum\limits_{n \in \mathbb{Z}} d_n e^{i n\omega},
\text{ where }
d_n:=c_{n}(Z_0).
\end{equation}
Let $k \in \mathbb{N}^{*}$ that will be chosen later on.
On the time interval $[0,k)$ we take
$w=0$, thus
$$Z(k^-,\omega)= Z_0(\omega) e^{i k \omega} =
\sum\limits_{n \in \mathbb{Z}} d_n e^{i(n+k)\omega}
\text{ and }
z(k^-,\omega)=z_0(\omega).$$

At time $k$, we apply the control $w(t)= i \pi \delta_{k}(t)$
in order to move close to $+e_3$. Indeed,
thanks to Definition \ref{defDirac-Cauchy-Bloch},
we have
$$M(k^+,\omega)=\exp(\pi \Omega_x)M(k^-,\omega)=
\left( \begin{array}{ccc}
1 & 0 & 0 \\ 0 & -1 & 0 \\ 0 & 0 & -1
\end{array} \right) M(k^-,\omega),$$
thus
$$Z(k^+,\omega)=\overline{Z(k^-,\omega)} =
\sum\limits_{n \in \mathbb{Z}} \overline{d_n} e^{i(-n-k)\omega}
\text{ and }
z(k^{+},\omega)=-z(k^-,\omega).$$

On the time interval $(k,3k)$ we apply a control of the form
$$w(t)=\sum\limits_{p=1}^{2k-1} w_{p} \delta_{p+k}(t),$$
where $w_p \in \mathbb{C}$. Approaching the nonlinear system
(\ref{eq:Z-sym}) by its linearized system around $(Z \equiv 0, w \equiv 0)$,
we get
\begin{equation} \label{approx-linearise}
\begin{array}{ll}
Z(3k^-,\omega)
& \sim
\Big( Z(k^+,\omega) - \int_{k}^{3k} w(t) e^{-i\omega(t-k)} dt \Big) e^{i 2k\omega}
\\ & \sim
\left( \sum\limits_{n \in \mathbb{Z}} \overline{d_{n}} e^{i(-n-k)\omega}
- \sum\limits_{p=1}^{2k-1} w_p e^{-i p \omega} \right) e^{i 2 k \omega}.
\end{array}
\end{equation}
Moreover, $z$ stays close to $+1$ because the control applied is small.
Choosing $w_p:=\overline{d_{p-k}}$, we get
$$Z(3k^-,\omega) \sim \sum\limits_{|n| \geqslant k} \overline{d_{n}} e^{i(-n+k)\omega}.$$

 Finally, at time $3k$, we apply the control $w(t)=i \pi \delta_{3k}(t)$
in order to return to $-e_3$:
$$Z(3k^+,\omega)=\overline{Z(3k^-,\omega)}
\sim \sum\limits_{|n| \geqslant k} d_n e^{i(n-k)\omega}$$
and $z(3k^+,\omega)=-z(3k^-,\omega)$ is close to $-1$.
Now, by choosing $k=k(Z_0)$ such that
$$\sum\limits_{|n| \geqslant k} |d_n| < \frac{1}{2} N(Z_0),$$
we get the existence of a time $T=T(Z_0):=3k$ and
a control $w:[0,T] \rightarrow \mathbb{C}$ such that
$N[Z(T)]<N[Z_0]/2$.
\\

 Finally, the steps that need to be justified are
\begin{itemize}
\item the approximation of the nonlinear system by its
linearized system in (\ref{approx-linearise}),
\item the convergence, for the norm $N$, of the solutions of (\ref{Bloch-ABC})
when we approximate the Dirac controls by controls in $L^{\infty}_{loc}$.
\end{itemize}

\subsection{Proof of Theorem \ref{Main-thm-(0,pi)}}
\label{subsec:proof of thm}

Let us recall that the solutions of (\ref{Bloch-ABC})
with Dirac controls have been defined in Definition \ref{defDirac-Cauchy-Bloch},
and that we have the following result.

\begin{Prop} \label{CV-rect-Dirac}
Let $\beta, \gamma \in \mathbb{R}$,
$M_0 \in C^0 ([0,\pi],\mathbb{S}^2)$ be such that $N(M_0)<+\infty$.
Let  $M$ be the solution of (\ref{Bloch-ABC}) with
$M(0)=M_0$, $u(t)=\beta \delta_0(t)$, and $v(t)= \gamma \delta_0(t)$.
 For $\epsilon>0$, let $M_{\epsilon}$ be the (classical) solution of
(\ref{Bloch-ABC}) with
$M(0)=M_0$, $u(t)=(\beta/\epsilon) 1_{(0,\epsilon)}(t)$, and
$v(t)= (\gamma/\epsilon) 1_{(0,\epsilon)}(t)$.
Then
\begin{gather}
\label{justiDirac}
N(M_{\epsilon}(\epsilon) - M(0^+) ) \rightarrow 0
\text{  when  } \epsilon \rightarrow 0.
\end{gather}
\end{Prop}

\textbf{Proof of Proposition \ref{CV-rect-Dirac}:}
We have
$$
\begin{array}{ll}
M_{\epsilon}(\epsilon,\omega)=
\exp[ \epsilon |\omega| \Omega_z + \beta \Omega_x + \gamma \Omega_y ] M_{0}(\omega),
&
M(0^+,\omega)=\exp[ \beta \Omega_x + \gamma \Omega_y ] M_{0}(\omega).
\end{array}
$$
One has
\begin{gather}
\label{estNMepsilon}
N(M_{\epsilon}(\epsilon)-M(0^+))\leqslant \sum_{n=1}^{+\infty}\frac{a_n(\epsilon)}{n!},
\end{gather}
with
\begin{gather}
\label{an=}
a_n(\epsilon):=N\Big(
[ (\epsilon |\omega|\Omega_z+\beta \Omega_x +\gamma \omega_y)^n
-(\beta \Omega_x +\gamma \omega_y)^n ] M_0(\omega)\Big).
\end{gather}
Noticing that $N(|\omega|)<+\infty$, using (\ref{N(fg)}) together with standard estimates and the Weierstrass $M$-test, one easily sees that (\ref{justiDirac}) follows from \eqref{estNMepsilon} and \eqref{an=}. $\Box$
\\

Thanks to Proposition \ref{CV-rect-Dirac},
Theorem \ref{Main-thm-(0,pi)} is a consequence of the following result.

\begin{thm} \label{Main-thm-precise}
There exists $\delta>0$ such that,
for every $M_0 : [0,\pi] \rightarrow \mathbb{S}^2$ with
$N[ Z_0 ] < \delta$ and $z_0 < -1/2$,
the solution of (\ref{Bloch-ABC}) with $M(0)=M_0$,
$$u(t):= \pi \delta_{k}(t) -
\sum\limits_{p=1}^{2k-1} \Im \Big(  c_{-k+p} (Z_0)  \Big)  \delta_{k+p}(t)  +
\pi \delta_{3k}(t)  ,$$
$$v(t):= - \sum\limits_{p=1}^{2k-1}
\Re \Big( c_{-k+p} (Z_0) \Big) \delta_{k+p}(t),$$
where $k=k(M_0) \in \mathbb{N}$ is such that (\ref{cdt-k}) holds,
satisfies
\begin{equation} \label{N[Z(3k+)]}
N[ Z(3k^+) ] < \frac{1}{2} N(Z_0),
\end{equation}
\begin{equation} \label{z(3k+)bis}
z(3k^+) < -\frac{1}{2}.
\end{equation}
\end{thm}

The key point of the proof of Theorem \ref{Main-thm-precise} is the
following result.

\begin{Prop} \label{Prop:1saut}
There exist $\mathcal{C}>0$ and  $\mathcal{C}'>0$ such that,
for every $d_0 \in \mathbb{C}$ with $|d_0| \leqslant 1$,
for every $M_0=(x_0,y_0,z_0):[0,\pi] \rightarrow \mathbb{S}^2$
with $N(Z_0) \leqslant 1$ and $z_0>0$,
the solution of (\ref{Bloch-ABC}) with $M(0)=M_0$,
$v(t)=-\Re( d_0 ) \delta_0(t)$, $u(t)= \Im ( d_0 ) \delta_0(t)$
satisfies
\begin{equation} \label{N2-0}
N \Big( Z(0^+) -   Z_0 + d_0  \Big)
\leqslant \mathcal{C} |d_0| \max \{ |d_0| , N(Z_0) \},
\end{equation}
\begin{equation} \label{z>}
z(0^+,\omega) \geqslant z_0(\omega) - \mathcal{C}'|d_0| \max\{ |d_0| , N(Z_0)\}.
\end{equation}
\end{Prop}

\textbf{Proof of Proposition \ref{Prop:1saut}:}
Let us write $d_0 = \beta_0 + i \gamma_0$, with $\beta_0, \gamma_0 \in \mathbb{R}$.
We have
$$M(0^+,\omega)=\exp[ \beta_0 \Omega_x + \gamma_0 \Omega_y] M_0(\omega).$$
Using the decomposition
\begin{equation} \label{dec-exp}
\exp [ \beta_0 \Omega_x + \gamma_0 \Omega_y ] = I + \beta_0 \Omega_x + \gamma_0 \Omega_y + R,
\text{ where } \| R \|= O( |d_0|^2 ) \text{ as } d_0\rightarrow 0,
\end{equation}
we get
$$Z(0^+,\omega)=Z_0(\omega) - d_0 z_0(\omega) +
R_1 x_0(\omega) + R_2 y_0(\omega) + R_3 z_0(\omega),$$
where $R_j \in \mathbb{C}$, $|R_j| \leqslant C |d_0|^2$ for $j=1,2,3$,
and $C$ is a universal constant. Therefore, we have
$$c_n[ Z(0^+) - Z_0 + d_0] =
d_0 c_n[ 1 - z_0] + R_1 c_n[x_0] + R_2 c_n[y_0] + R_3 c_n[z_0], \forall n \in \mathbb{Z}.$$
Using (\ref{NxyZ}) and (\ref{N(z-1)}), we get
$$\begin{array}{ll}
N \Big( Z(0^+) -  Z_0 + d_0  \Big)
& \leqslant
|d_0| N(z_0-1)+ |R_1| N(x_0) + |R_2| N(y_0) + |R_3| N(z_0)
\\ & \leqslant
2 |d_0| N(Z_0)^2 + C |d_0|^2 [ 2 N(Z_0) + 1 + 2 N(Z_0)^2 ],
\end{array}$$

which gives (\ref{N2-0}) with $\mathcal{C}=2+5C$.

 From (\ref{dec-exp}) we get
$$z(0^+,\omega)=z_0(\omega)+ \Re \Big( \overline{d}_0 Z_0(\omega) \Big)
+ R_1' x_0(\omega) + R_2' y_0(\omega) + R_3' z_0(\omega),$$
where $R_j' \in \mathbb{C}$, $|R_j'| \leqslant C' |d_0|^2$ for $j=1,2,3$,
where $C'$ is another universal constant.
Using (\ref{NxyZ}) and (\ref{N(z-1)}), we get
$$\begin{array}{ll}
z(0^+,\omega)
& \geqslant
z_0(\omega) - |d_0| |Z_0(\omega)| -
C' |d_0|^2 [ |x_0(\omega)| + |y_0(\omega)| + |z_0(\omega)| ]
\\ & \geqslant
z_0(\omega) - |d_0| N(Z_0) -
C' |d_0|^2 [ N(x_0) + N(y_0) + N(z_0) ]
\\ & \geqslant
z_0(\omega) - |d_0| N(Z_0) - C' |d_0|^2 [ 2 N(Z_0) + 1 + 2 N(Z_0)^2]
\end{array}
$$
which gives (\ref{z>}) with $\mathcal{C}'=1+5C'$.$\Box$
\\

\textbf{Proof of Theorem \ref{Main-thm-precise}:}
Let $\delta$ be such that
\begin{equation} \label{hyp-delta-C}
\begin{array}{ccc}
4 \mathcal{C} \delta < 1,
&
\mathcal{C}' \delta < 1/2,
&
\delta \in (0,1),
\end{array}
\end{equation}
where
$\mathcal{C}, \mathcal{C}'$ are as in Proposition \ref{Prop:1saut}.
Let $M_0$, $k$, $u$, $v$ be as in Theorem \ref{Main-thm-precise}.
We use the notation (\ref{dec:Z0}).
\\

\emph{First step: on [0,k].}
We have (see the previous section)
\begin{equation} \label{Zz(k+)}
Z(k^+,\omega)=
\sum\limits_{n \in \mathbb{Z}} \overline{d_n} e^{i(-n-k)\omega}
\text{  and  }
z(k^+,\omega)=-z_0(\omega).
\end{equation}
\\

\emph{Second step: on (k,3k).}
Let us prove by induction on $p \in \{0,...,2k-1\}$ that
for every $p \in \{0,...,2k-1\}$, we have
\begin{equation}
\begin{array}{ll}
(H_p) \text{: }
&
N \Big[
Z \Big( (k+p)^+ \Big) -
\sum_{n \in \mathbb{Z}-\{-k+1,...,-k+p \}}
\overline{d_n} e^{i(-n-k+p)|\omega|}
\Big]
\\ & \leqslant
\mathcal{C} [ |d_{-k+1}| + ... + |d_{-k+p|} ] N(Z_0),
\end{array}
\end{equation}
\begin{equation}
(H'_p) \text{:  }
z((k+p)^+,\omega) \geqslant -z_0(\omega)
- \mathcal{C}' [ |d_{-k+1}|+ ... + |d_{-k+p}|]N(Z_0).
\end{equation}

Notice that $(H_{2k-1})$ and (\ref{hyp-delta-C}) provide
$$N \Big[ Z \Big( 3k^- \Big) -
\sum_{n \in \mathbb{Z}, |n| \geqslant k} \overline{d_n} e^{i(-n+k)|\omega|}
\Big] \leqslant \mathcal{C} N(Z_0)^2,$$
thus, thanks to (\ref{hyp-delta-C}) and (\ref{cdt-k}), we have
\begin{equation} \label{N[Z(3k-)]}
N[Z(3k^-)] < N[Z_0]/2.
\end{equation}
We also have, thanks to $(H'_{2k-1})$ and (\ref{hyp-delta-C})
$$z(3k^-,\omega)=z((3k-1)^+,\omega) \geqslant
-z_0(\omega) + \mathcal{C}' \delta^2 > 0$$
thus
\begin{equation} \label{z(3k-)}
z(3k^-,\omega) = \sqrt{1-|Z(3k^+,\omega)|^2}
> \sqrt{1-\delta/2} > 1/2.
\end{equation}

The properties $(H_0)$ and $(H'_0)$ come from (\ref{Zz(k+)}).
Now, let $p\in \{1,...,2k-1\}$ and let us assume that $(H_{p-1})$
and $(H'_{p-1})$ hold. Thanks to $(H_{p-1})$ and (\ref{hyp-delta-C}), we have
$$N[ Z((k+p)^-) ] = N[ Z((k+p-1)^+) ] \leqslant N[Z_0]
\leqslant \delta \leqslant 1$$
and thanks to $(H'_{p-1})$ we have
$$z((k+p)^-,\omega)=z((k+p-1)^+,\omega)\geqslant
-z_0(\omega) - \mathcal{C}' N(Z_0)^2 >
\frac{1}{2} - \frac{1}{2}=0$$
thus we can apply Proposition \ref{Prop:1saut}.
Thanks to Proposition \ref{Prop:1saut} and $(H_{p-1})$, we get
$$\begin{array}{ll}
& N \Big[  Z \Big( (k+p)^+ \Big) -
\sum\limits_{n \in \mathbb{Z}-\{-k+1,...,-k+p\}} \overline{d_n} e^{i(-n-k+p)\omega}
\Big]
\\ \leqslant &
N \Big[  Z \Big( (k+p)^+ \Big) -   Z \Big( (k+p)^- \Big) + d_{-k+p} \Big]
\\ &
+ N \Big[ Z \Big( (k+p)^- \Big) -
\sum\limits_{n \in \mathbb{Z}-\{-k+1,...,-k+p-1\}} \overline{d_n}   e^{i(-n-k+p)\omega}
\Big]
\\ \leqslant &
\mathcal{C} |d_{-k+p}| N[ Z((k+p)^-) ]
\\ &
+ N \Big[ Z \Big( (k+p-1)^+ \Big) -
\sum\limits_{n \in \mathbb{Z}-\{-k+1,...,-k+p-1\}} \overline{d_n}   e^{i(-n-k+p-1)\omega}
\Big]
\\ \leqslant &
\mathcal{C} |d_{-k+p}| N [ Z_0 ]
+ \mathcal{C} [ |d_{-k+1}| + ... + |d_{-k+p-1}| ] N(Z_0),
\end{array}$$
which proves $(H_{p})$. Thanks to Proposition \ref{Prop:1saut} and $(H'_{p-1})$, we get
$$\begin{array}{ll}
z((k+p)^+,\omega)
& \geqslant
z((k+p)^-,\omega) - \mathcal{C}' |d_{-k+p}| N[ Z((k+p)^-) ]
\\ & \geqslant
z((k+p-1)^+,\omega) - \mathcal{C}' |d_{-k+p}| N[ Z_0 ]
\\ & \geqslant
-z_0(\omega) - \mathcal{C}'[|d_{-k+1}+...+|d_{-k+p-1}|+|d_{-k+p}|] N[Z_0].
\end{array}$$

\emph{Third step: at 3k.} We have
$Z(3k^+,\omega)=\overline{Z(3k^-,\omega)}$ and
$z(3k^+,\omega)=-z(3k^-,\omega)$,
thus (\ref{N[Z(3k-)]}) and (\ref{z(3k-)})
give (\ref{N[Z(3k+)]}) and (\ref{z(3k+)bis}). $\Box$

\section{Comparison}
\label{sec:Comparaison}

In this section, we compare the control results and processes presented in sections
\ref{Li-Khaneja} and \ref{Section:asymptotic}.
\\

 First, let us compare the statements of Theorems \ref{thm:LK}
(or Corollary 1) and \ref{Main-thm-(0,pi)}.
On one hand, the statement of Theorem \ref{thm:LK} is stronger than
the one of Theorem \ref{Main-thm-(0,pi)} because it is global and it gives the
approximate controllability of (\ref{Bloch-ABC})
for the norms  $\|.\|_{H^s}, \forall s < 1$
(whereas Theorem \ref{Main-thm-(0,pi)} only provides the approximate controllability
for $N$). On the other hand, Theorem \ref{Main-thm-(0,pi)}
is stronger than Theorem \ref{thm:LK}
because it needs less regular initial data.
\\

Now, let us compare the control processes detailed in the proof of
Theorems \ref{thm:LK} and \ref{Main-thm-(0,pi)}.
Given $M_0 \in H^1((\omega_*,\omega^*),\mathbb{S}^2)$,
the proofs of Lemma \ref{Lem:dec} and Proposition \ref{Prop:LK} give an explicit way
to find $T>0, u,v, \in D$ such that
$$\| U[T^+;u,v,M_0]' \|_{L^2} < \| M_0' \|_{L^2}.$$
Iterating this process, we produce a sequence of reachable points
$(M_{n})_{n \in \mathbb{N}} \subset H^1((\omega_*,\omega^*),\mathbb{S}^2)$
such that $(\|M_n' \|_{L^2})_{n \in \mathbb{N}}$ decreases.
We expect that $\| M_n'\|_{L^2} \rightarrow 0$ when $n \rightarrow + \infty$,
and once this norm is small enough, we apply a control given in
Lemma \ref{Lem:dec} \textbf{(2)} to go closer to $e_3$.
However, the sequence $(M_n')_{n \in \mathbb{N}}$ may not converge to $0$.
Thus, the control process presented in section
\ref{Li-Khaneja} is not completely satisfying from a practical point of view.

Moreover, even if the sequence $(M_n')_{n \in \mathbb{N}}$ converges to
$0$ in $L^2((\omega_*,\omega^*),\R^3)$, the controllability process may
take a long time (in particular the controllability time is not a priori
bounded by a quantity depending only on $\|M_0\|_{H^1 }$)
and cost a lot (because at each step, one has to compute new controls $u$ and $v$
and because the commands proposed in the proof of Proposition \ref{Prop:LK} involve
many trips between $-e_3$ and $+e_3$).

On the contrary, the controllability process presented in section
\ref{Section:asymptotic} works within a time $T$ which is explicit, with controls $u$, $v$ that are also explicit in terms of
the Fourier coefficients of $M_0$, and needs only two trips between $\pm e_3$.
Thus, the time and the cost are well known.
\\

Let us compare the time and the cost involved by the two
controllability processes on a particular example.
We take $(\omega_*,\omega^*)=(0,\pi/2)$, and an initial data of the form
$$M_{0}(\omega)=
\left( \begin{array}{c}
\epsilon x_{\epsilon} (\omega) \\ 0 \\ \sqrt{1 - \epsilon^2 x_{\epsilon}(\omega)^2}
\end{array} \right)$$
where $\epsilon>0$ is small,
\begin{equation} \label{x-epsilon}
x_{\epsilon}(\omega) = \sum\limits_{k=1}^N a_k(\epsilon) \cos((2k-1) \omega)
+ \cos((2N+1)\omega), \forall \omega \in (0,\pi/2),
\end{equation}
and $(a_{k}(\epsilon))_{1 \leqslant k \leqslant N} \in \mathbb{R}^N$ are such that
\begin{equation} \label{og-pol}
\int_{0}^{\pi/2} \frac{x_{\epsilon}'(\omega)}{\sqrt{1-\epsilon^2 x_{\epsilon}(\omega)^2}}
\omega^K d\omega=0, \forall K \in \{0,...,N-1\}.
\end{equation}
We will prove later the existence of such coefficients. We want to reach $e_{3}$.
\\

Let us apply the strategy presented in section \ref{Li-Khaneja},
to find explicit $T>0$, $u,v \in D$ such that
$$\| U[T^+;u,v,M_0]' \|_{L^2} < \| M_0' \|_{L^2}.$$
One needs a polynomial $Q \in \mathbb{R}[X]$ such that
$$\int_0^{\pi/2} (zx'-z'x)Q d\omega <0.$$
Then, $\text{deg}(Q) \geqslant N$, because of (\ref{og-pol}).
Thanks to the proof of Lemma \ref{Lem:dec},
there exists $\tau^*=\tau^*(Q,x_{\epsilon})>0$ and $\alpha>0$ such that,
for every $\tau \in (0,\tau^*)$, there exist
$T >0$, $u,v \in D$ such that
$$\| U[T^+;u,v,M_0]' \|_{L^2}^2 \leqslant \| M_0' \|_{L^2}^2 - \alpha \tau.$$
However $\tau^*$ cannot be quantified, thus,
we do not know the size of the decrease.
Moreover, as emphasized in the proof of Proposition \ref{Prop:LK},
the time of control $T$ satisfies $T \geqslant 2^N \tau^{1/N}$
(time needed to generate $I+\tau \omega^N \Omega_x + o(\tau)$)
and one makes more than $2^N$ trips between $\pm e_3$
(just count how many times the matrices
$\exp( \pi \Omega_x)$ or $\exp( \pi \Omega_y)$ appear in
the generation of $I+\tau \omega^N \Omega_x + o(\tau)$ in the
proof of Proposition \ref{Prop:LK}).
\\

With the strategy of section \ref{Section:asymptotic}
taking the same explicit expression for $M_0$ on $(-\pi,\pi)$,
we know the existence of $\epsilon^*>0$ such that,
for every $\epsilon \in (0,\epsilon^*)$,
the explicit controls
$$u(t):= \pi \delta_{2N+1}(t) + \pi \delta_{6N+3}(t),$$
$$v(t):=
-\frac{\epsilon}{2} \sum_{m=1}^{N} a_{N+1-m}(\epsilon) \delta_{2N+1+2m}(t)
-\frac{\epsilon}{2} \sum_{m=N+1}^{2N} a_{m-N}(\epsilon) \delta_{2N+1+2m}(t),$$
with the convention $a_{N+1}(\epsilon)=a_{-N-1}(\epsilon)=1$, realize
$$N \Big[ U[(6N+3)^+;u,v,M_{0}] + e_3 \Big] <
\frac{1}{2} N[M_0+e_3].$$
Here, the controls are explicit,
the time scales like $6N$,
we have a bound from below for the decrease of the $N$-distance to $e_3$,
and the process needs only 2 trips between $\pm e_3$.
\\

Now, let us prove the existence of the coefficients
$(a_{k}(\epsilon))_{1 \leqslant k \leqslant N}$.

\begin{Lem}
Let $N \in \mathbb{N}^*$.

\textbf{(i)} The matrix $A \in \mathcal{M}_N(\mathbb{R})$ with coefficients
$$A_{k,K} := \int_0^{\pi/2} (2k-1) \sin((2k-1)\omega) \omega^K d\omega,
1 \leqslant k \leqslant N, 0 \leqslant K \leqslant N-1,$$
is invertible.

\textbf{(ii)} There exists $\epsilon^*>0$ and a $C^1$ map
$\epsilon \in [0,\epsilon^*] \mapsto
(a_k(\epsilon))_{1 \leqslant k \leqslant N} \in \mathbb{R}^N$
such that (\ref{x-epsilon})-(\ref{og-pol}) hold.
\end{Lem}

\textbf{Proof:} \textbf{(i)} We assume that $A$ is not invertible.
Then, there exists $(\lambda_1,...,\lambda_N) \in \mathbb{R}^N-\{0\}$ such that
\begin{equation} \label{hyp-absurd}
\int_0^{\pi/2} \sum_{k=1}^N \lambda_k \sin((2k-1)\omega) \omega^K d\omega =0,
\forall 0 \leqslant K \leqslant N-1.
\end{equation}
Let $f(\omega):=\sum_{k=1}^N \lambda_k \sin((2k-1)\omega)$
and $0<\omega_1<...<\omega_L<\pi/2$ be all the values of
the open interval $(0,\pi/2)$ on which
$f$ vanishes and changes its sign. Then, the function
$\omega \mapsto f(\omega)(\omega-\omega_1)...(\omega-\omega_L)$
has a constant sign on $(0,\pi/2)$ and it is not identically zero, thus
$$\int_0^{\pi/2} f(\omega)(\omega-\omega_1)...(\omega-\omega_L) d\omega \neq 0.$$
The assumption (\ref{hyp-absurd}) ensures that $L \geqslant N$.
Thanks to trigonometric formulas, there exists
$(\mu_1,...,\mu_N) \in \mathbb{R}^N-\{0\}$ such that
$$f(\omega)=\sum_{k=1}^{N} \mu_k \sin(\omega)^{2k-1}=
\sin(\omega) \sum_{k=1}^{N} \mu_k \sin(\omega)^{2(k-1)}.$$
Since the quantities $\sin(\omega_1)^2,...,\sin(\omega_N)^2$
are all different from zero ($\omega_1,...\omega_N \in (0,\pi/2)$),
they provide $N$ roots for the polynomial
$$\sum_{k=1}^{N} \mu_k X^{(k-1)}$$
that have a degree $\leqslant (N-1)$ and is different from zero.
This is a contradiction.

\textbf{(ii)} Thanks to \textbf{(i)}, there exists
$(\alpha_1,...,\alpha_N) \in \mathbb{R}^N$ such that
$$\int_0^{\pi/2} \left(
\sum_{k=1}^{N} \alpha_k (2k-1) \sin((2k-1)\omega)
+ (2N+1) \sin((2N+1)\omega) \right) \omega^K d\omega =0,
\forall 0 \leqslant K \leqslant N.$$
There exists $M>0$ such that
$$\Big| \sum_{k=1}^{N} \alpha_k \cos((2k-1)\omega) + \cos((2N+1)\omega) \Big| \leqslant M,
\forall \omega \in (0,\pi/2).$$
When $b=(b_1,...,b_N)^t \in \mathbb{R}^N$ we have
$$\Big| \sum_{k=1}^N (\alpha_k+b_k) \cos((2k-1)\omega) + \cos((2N+1)\omega) \Big|
\leqslant M + \sqrt{N} \|b\|$$
thus, the following map $F$ is well defined
$$\begin{array}{crclccc}
 F:& \left( 0 , \frac{1}{2M} \right)  & \times&
B_{\mathbb{R}^N}\left( 0 , \frac{M}{\sqrt{N}} \right)  & \rightarrow & \mathbb{R}^N \\
  & (\epsilon  & ,     &  b),                             & \mapsto     & F(\epsilon,b)
\end{array}$$
$$F(\epsilon,b):=\Big(
\int_0^{\pi/2} \frac{y_b'(\omega)}{\sqrt{1-\epsilon^2 y_b(\omega)^2}} \omega^K d\omega
\Big)_{1 \leqslant K \leqslant N}$$
where
$$y_b(\omega):=\sum_{k=1}^N (\alpha_k+b_k)\cos((2k-1)\omega) +
\cos((2N+1)\omega), \forall \omega \in (0,\pi/2).$$
Then $F(0,0)=0$ and $d_b F(0,0)$ is invertible, thanks to \textbf{(i)}.
Thus, the implicit function theorem gives the conclusion. $\Box$

\end{document}